%Paper: alg-geom/9311013
%From: Takao Fujita <fujita@math.titech.ac.jp>
%Date: Tue, 30 Nov 93 17:33:12 JST

%amstex version 1.1c
\input amstex
\documentstyle{amsppt}
\nologo
\magnification=\magstep1
\hsize=145mm
\vsize=220mm
\hcorrection{-7mm}
\vcorrection{-10mm}
\baselineskip=18pt
 \abovedisplayskip=4pt
 \belowdisplayskip=4pt
\parskip=3pt
\parindent=8mm
%\input amsfonts

%The above may be rewritten according to your taste.

\topmatter
\title Remarks on Ein-Lazarsfeld criterion \\
of spannedness of adjoint bundles of polarized threefold
\endtitle
\author Takao FUJITA \endauthor
\address {Takao FUJITA\newline
Department of Mathematics
\newline
Tokyo Institute of Technology
\newline
Oh-okayama, Meguro, Tokyo
\newline
152 Japan
\newline
e-mail:fujita\@math.titech.ac.jp}
\endaddress
\endtopmatter

\define\almb{\allowmathbreak }
\define\dnl{\newline\newline}
\define\lra{{\longrightarrow}}
\define\mb{\mathbreak }
\define\nl{\newline}
\define\nomb{\nomathbreak }

\define\Bs{{\roman{Bs}}}

\define\Inf{{\roman{Inf}}}
\define\Int{{\roman{Int}}}

\define\Max{{\roman{Max}}}
\define\Min{{\roman{Min}}}
\define\Pic{{\roman{Pic}}}

\define\dm{{\roman{dim}}}
\define\rk{{\roman{rank}}}

\define\SO{\Cal O}
\define\SP{\Cal P}

\define\BC{{\Bbb C}}
\define\BP{{\Bbb P}}
\define\BQ{{\Bbb Q}}

\define\Cor{{\bf Corollary. }}

\define\Lem{{\bf Lemma. }}
\define\Prf{{\it Proof. }}

\define\Rmk{{\it Remark. }}
\define\Th{{\bf Theorem. }}

\document
\noindent{\bf Introduction}

Let $L$ be an ample line bundle on a smooth complex projective variety $M$ of
dimension $n$.
We conjecture that $K+nL$ is spanned unless $L^n=1$, where $K$ is the canonical
bundle of $M$.
In this paper we prove this conjecture for $n=3$.

For $n=2$, the conjecture follows easily from Reider's theory [{\bf R}].
In higher dimensions, after the notable works [{\bf Kol}] and [{\bf Dem}], a
breakthrough has been achieved by [{\bf EL}] in case $n=3$.
They proved, among others, the following criterion:

{\sl Let $B$ be a nef and big line bundle on a smooth 3-fold $M$, and let $x$
be a point on $M$.
Suppose that\nl
$BC\ge3$ for any curve $C$ in $M$ with $C\ni x$,\nl
$B^2S\ge7$ for any irreducible surface $S$ in $M$ with $S\ni x$, and\nl
$B^3\ge92$.\nl
Then $x\notin\Bs{\vert}K+B{\vert}$, i.e., $K+B$ is spanned at $x$.}

The spannedness of $K+3L$ follows from this when $L^4\ge4$.
Moreover, with some additional techniques, they proved that $K+4L$ is always
spanned.

The purpose of this paper is to show that the above criterion holds true even
if the last assumtion $B^3\ge92$ is replaced by $B^3\ge51$.
For the proof we use the same techniques as in [{\bf EL}].
Our argument starts from the following

{\bf Theorem 0.} {\sl Let $M$, $B$, $x$ be as above, let $\pi_1:M_1\lra M$ be
the blow-up at $x$, let $E$ be the exceptional divisor of it and let
$\sigma_1$, $\sigma_2$, $\sigma_3$ be positive rational numbers such that
$$\sigma_2\ge\Max(\frac{2\sigma_3}{\sigma_3-1},\frac{\sqrt2\sigma_3}{\sigma_3-2})\text{,\quad }
\sigma_1\ge\Max(\frac{2\sigma_3}{\sigma_3-1},
\frac{\sigma_3}{\sigma_3-3})\text{,\quad}
\sigma_3>3.$$
Suppose that $BC>\sigma_1$ for any curve $C\ni x$, $B^2S>\sigma_2^2$ for any
surface $S\ni x$ and that $\pi_1^*B-\sigma_3E$ is big.
Then $K+B$ is spanned at $x$.}

In their Theorem $1^*$ in [{\bf EL}], they assumed $B^3>\sigma_3^3$ instead of
the bigness of $\pi_1^*B-\sigma_3E$.
The bigness follows from this assumption, but thereafter this is never used.
Thus, [{\bf EL}] proves our Theorem 0 implicitly.

We will apply Theorem 0 for $\sigma_1=3$, $\sigma_2=\sqrt 7$ and
$\sigma_3=\frac92$.
It suffices to show that the bigness follows from $B^3\ge51$, if $K+B$ is not
spanned at $x$.
To this end, we need an estimate of the rank of the restriction map $H^0(M_1,
s\pi_1^*B-jE)\lra H^0(E,\SO_E(j))$.
Of course $\le h^0(\BP^2,\SO(j))=(j+1)(j+2)/2$, which is enough when
$B^3>\sigma_3^3=91.125$.
However, when $K+B$ is not spanned at $x$, we get a sharper bound by a careful
quantitative analysis, which yields the bigness of $\pi_1^*B-\frac92E$ in case
$B^3\ge51$.

This paper is organized as follows.
In \S1, we describe the basic situation and fix notation.
\S2 and \S3 are devoted to the method for giving the above mentioned estimate.
In \S4 we complete the proof of main results and discuss related topics.

Needless to say, I benefited very much by reading the preprints [{\bf EL}],
[{\bf Dem}], [{\bf Kol}], [{\bf T}].
I would like to express my hearty thanks to the authors, especially to
Professors Ein and Lazarsfeld.
\dnl
{\bf \S1. Basic observations}

(1.1) Usually we follow the customary notation in modern algebraic geometry,
and that in [{\bf EL}].

A line bundle $L$ on a variety $V$ is said to be {\it big} if $\kappa(L)=n=\dm
V$, i.e., $h^0(V,tL)$ is a function in $t$ of the growth order $t^n$ when
$t\lra\infty$.

The integral part of a $\BQ$-divisor $D$ will be denoted by $\Int(D)$.
Thus $D-\Int(D)$ is the fractional part.

For $i=1,2$, let $D_i$ be a prime Weil divisor on a normal variety $V_i$ which
is birational to $V$.
Such pairs $(V_i, D_i)$ determine the same discrete valuation of the function
field of $V$ if and only if there exist another model $\tilde V$ with
birational morphisms $f_i:\tilde V\lra V_i$ and a prime divisor $\tilde D$ on
$\tilde V$ such that $f_i(\tilde D)=D_i$.
Letting $(V_1,D_1)\sim(V_2,D_2)$ in such cases, $\sim$ is an equivalence
relation.
An equivalence class with respect to $\sim$ will be called a {\it place} of
$V$.
The place represented by such a pair $(V,D)$ will be denoted by $\SP(V,D)$ or
simply by $\SP(D)$.

For a place $P$ of $V$, take a representative pair $(V',D')$ such that
$\pi:V'\lra V$ is a morphism.
The image $Z=\pi(D')$ is independent of the choice of the pair, so it will be
called the {\it locus} of $P$ on $V$.
We say also that $P$ lies over the locus $Z$.

For a subvariety $Z$ of $V$ not contained in the singular locus of $V$, let
$V'$ be the blowing-up of $V$ along $Z$ and let $E_Z$ be the exceptional
divisor over $Z$.
The place $\SP(E_Z)$ will be called the primary place over $Z$.

(1.2) Let $P=\SP(V',D')$ be a place of a variety $V$ as above.
For a Weil divisor $D$, let $v_P(D)$ denote the coefficient of $D'$ in $f^*D$.
For $\psi\in H^0(V,L)$, we set $v_P(\psi)=v_P(D_\psi)$, where $D_\psi$ is the
member of ${\vert}L{\vert}$ corresponding to $\psi$.
Here we set $v_P(\psi)=+\infty$ by convention when $\psi=0$.

A linear system $\Lambda$ on $V$ can be identified with a linear subspace
$V(\Lambda)$ of $H^0(V,L)$.
We set $v_P(\Lambda)=\Min_{\psi\in V(\Lambda)}\{v_P(\psi)\}$.
For a place $Q$ of a subvariety $W$ of $V$, we set
$v_Q(\Lambda)=v_Q(\Lambda_W)$, where $\Lambda_W$ is the linear system
corresponding to the image of $V(\Lambda)$ via the restriction map
$H^0(V,L)\lra H^0(W,L_W)$.
Similarly we define $v_Q(\psi)=v_Q(\psi_W)$ and $v_Q(D)=v_Q(D_W)$.
In particular, for an effective divisor $D$, $v_Q(D)=+\infty$ if and only if
$W\subset D$.

We set $\Xi_Q(L)=\Inf_{t>0}\{v_Q({\vert}tL{\vert})/t\}$, which will be called
the {\it algebraic Lelong number} of $L$ at $Q$.
We have $\Xi_Q(L_W)\le\Xi_Q(L)$, but the equality does not hold in general,
since $H^0(V,tL)\lra H^0(W,tL_W)$ is not always surjective.
When $Q$ is the primary place over $Z$, $\Xi_Q(L)$ is sometimes denoted by
$\Xi_Z(L)$.

Clearly we have $\Xi_Q(mL)=m\Xi_Q(L)$ for any positive integer $m$.
Hence, for any $\BQ$-bundle $F$, $\Xi_Q(F)$ is well-defined in a natural way.
Since $\Xi_Q(m_1L_1+m_2L_2)\le m_1\Xi_Q(L_1)+m_2\Xi_Q(L_2)$ for any positive
integers $m_1, m_2$, $\Xi_Q(F)$ is a convex function in
$F\in\Pic(V)\otimes\BQ$.

(1.3) Now we consider a situation as in the introduction.
$B$ is a nef and big line bundle on a smooth 3-fold $M$, $x$ is a point on $M$
and $BC\ge3$ for any curve $C\ni x$, $B^2S\ge7$ for any surface $S\ni x$ and
$d=B^3\ge51$.
Let $\pi_1:M_1\lra M$ be the blow-up at $x$ and let $E$ be the exceptional
divisor over $x$.
Then $h^0(M_1,s\pi_1^*B-tE)-h^0(M_1,s\pi_1^*B-(t+1)E)\le
h^0(E,\SO(t))=(t+1)(t+2)/2$, so
$h^0(M_1,s\pi_1^*B-smE)\ge
h^0(M_1,s\pi_1^*B)-\sum_{t=0}^{sm-1}(t+1)(t+2)/2=h^0(M,sB)-sm(sm+1)(sm+2)/6$.
Since $h^0(M,sB)$ grows like $ds^3/6$ when $s\lra\infty$, we infer that
$\pi_1^*B-mE$ is big if $d>m^3$, in particular if $m=3$.

Now we have:

(1.4) {\sl There is a positive integer $\tau$ such that
$v_E({\vert}\tau(\pi_1^*B-3E){\vert})=0$}.

Indeed, since ${\vert}a(\pi_1^*B-3E){\vert}\ne\emptyset$ for some $a>0$, we
have $v_E({\vert}a\pi_1^*B-wE{\vert})=0$ for some $w\ge3a$.
On the other hand, $B$ is almost base point free in the sense of Goodman (cf.
[{\bf F};(6.13)]), so $v_E({\vert}\ell\pi_1^*B-\delta E{\vert})=0$ for some
$\ell, \delta>0$ with $3\ell>\delta$.
Then $v_E({\vert}L{\vert})=0$ for
$L=(3\ell-\delta)(a\pi_1^*B-wE)+(w-3a)(\ell\pi_1^*B-\delta
E)=(w\ell-a\delta)(\pi_1^*B-3E)$, so $\tau=w\ell-a\delta$ has the desired
property.

Take a sufficiently large integer $\tau $ as above.
Then we have a birational morphism $\tilde\pi:\tilde M\lra M_1$ satisfying the
following conditions:\nl
(1) The moving part $H$ of $\Lambda=\tilde\pi^*{\vert}\tau
(\pi_1^*B-3E){\vert}$ has no base point and $H^3>0$.\nl
(2) Let $F$ be the fixed part of $\Lambda$ and let $R$ be the ramification
divisor of $\pi=\pi_1\circ\tilde\pi:\tilde M\lra M$.
Then $F+R$ is supported on a simple normal crossing divisor.

Such a decomposition $\Lambda=F+{\vert}H{\vert}$ will be called a Hironaka
decomposition (or model) of degree $\tau $.
There are many such decompositions of various degrees.
The choice of the degree will be made precise later, so that some additional
conditions are satisfied.

(1.5) In the above situation, let $R=\sum_ja_jF_j$ and $\tau
^{-1}F+3\tilde\pi^*E=\sum_j\nu_jF_j$, where $F_j$'s are prime components of
$F+R$.
Sometimes $a_j$ (resp. $\nu_j$) will be denoted by $a(F_j)$ (resp. $\nu(F_j)$).
Put $c(F_j)=c_j=(a_j+1)/\nu_j$ for each $j$.

Let $F_0$ be the proper transform of $E$.
Then $\nu_0=\tau ^{-1}v_{F_0}(F)+3v_{F_0}(\tilde\pi^*E)=3$ by (1.4) and
$a_0=2$,
hence $c_0=1$.

Suppose that there is a component $F_j$ such that $c_j\le1$ and $F_j\cap
F_0\ne\emptyset$.
The place of $E$ represented by $(F_0, F_0\cap F_j)$ will be called a {\it bad
place} from now on.

(1.6) If there is no bad place as above, set $\sum_j
(\nu_j-a_j)F_j=P-N+\Delta$, where $\Delta$ is the fractional part of the left
hand side and $P$, $N$ are effective divisors without common components.
By Kawamata-Viehweg's vanishing theorem $H^1(\tilde M,\pi^*(K+B)-P+N)=0$, since
$\pi^*(K+B)-P+N-\tilde K-\Delta=\tau ^{-1}H$ is nef and big.
Since $\nu_0-a_0=1$, $P'=P-Z$ is effective.
Moreover $P'$ does not meet $F_0$ by the assumption.

On the other hand, any component $F_j$ of $N$ is $\pi$-exceptional since
$a_j>0$.
Hence, argueing similarly as in the case [{\bf EL};(4.1)], we infer that the
restriction map $H^0(\tilde M,\pi^*(K+B)-P')\almb\lra
H^0(F_0,\pi^*(K+B)_{F_0}-P'_{F_0})\cong\BC$ is surjective.

This implies $x\notin\Bs{\vert}K+B{\vert}$ and we are done in this case.

(1.7) When there is a bad place $P=\SP(F_0\cap F_j)$, set
$\nu'_j=v_{F_j}(\tilde\pi^*E)$ and $\nu''_j=v_{F_j}(\tau ^{-1}F)$.
Then $\nu_j=3\nu'_j+\nu''_j$.
Hence $1\ge c_i$ implies $a_i+1\le\nu_i=3\nu'_i+\nu''_i$.

Let $R_0$ be the ramification divisor of the map $F_0\lra E$.
Then $R_0=\sum_j a_jF_j+F_0-3\tilde\pi^*E$ in $\Pic(F_0)$,
so $v_P(R_0)=a_i-3\nu'_i\le\nu''_i-1$.
We divide the cases according to the nature of the place $P$.
\dnl
{\bf \S2. The case of non-exceptional place}

(2.1) Let things be as in (1.7) and suppose that the locus of $P$ on
$E\cong\BP^2$ is a curve.
In this case $v_P(R_0)=0$, so $\nu''_i\ge 1$.

Assume that $\Xi_P(\pi_1^*B-3E)<1$.
Then, by the definition of algebraic Lelong number, there is a multiple $\tau'$
of $\tau $ such that $v_P(F')/\tau'<1$ for the fixed part $F'$ of
$\tilde\pi'{}^*{\vert}\tau'(\pi_1^*B-3E){\vert}$ on $\tilde M'$.
Replacing the Hironaka decomposition by another one of degree $\tau'$, we
obtain $\nu''_i<1$, and we can get rid of such a situation.
Therefore we assume $\Xi_P(\pi_1^*B-3E)\ge1$ here.

(2.2) For each rational number $x>0$, put
$\xi(x)=\Xi_P(\pi_1^*B-xE)$; here put $\xi(x)=+\infty$ if
$\kappa(\pi_1^*B-xE)<0$.
This is a convex function in $x$, so
$\xi(x)\ge\eta_\lambda(x)=_{\text{def}}\Max(\lambda(x-3)+1,0)$ for some
$\lambda\ge\frac13$.

Let us now consider the restriction map
$\rho_{s,j}: H^0(M_1,s\pi_1^*B-jE)\lra H^0(E,\SO_E(j))$.
For any $\psi\in{\roman{Im}}(\rho_{s,j})$, we have
$v_P(\psi)\ge s\xi(\frac js)$
by the definition of $\xi$.
Hence $\rk(\rho_{s,j})\le\mb h^0(E,\SO(j-s\xi(\frac
js)))\le(1+j-s\eta_\lambda(\frac js))(2+j-s\eta_\lambda(\frac js))/2$,
provided $1+j-s\eta_\lambda(\frac js)\ge 0$.

(2.3) Before proceeding further, we set\nl
$d_1(s,m,\lambda)=\sum_{j=0}^{sm-1}\frac 1s=m$,\nl
$d_2(s,m,\lambda)={\dsize\sum}_{j=0}^{sm-1}\dfrac32(\dfrac js
-\eta_\lambda(\dfrac js))\dfrac 1s$ and\nl
$d_3(s,m,\lambda)={\dsize\sum}_{j=0}^{sm-1}\dfrac12(\dfrac js
-\eta_\lambda(\dfrac js))^2\dfrac1s$.\nl
Here, $sm$ is not always an integer, under the following convention:
Given a sequence $x_0$, $x_1$, $\cdots$ and a real number $r\ge0$, we set
$\sum_{j=0}^r=_{\text{def}}\sum_{j=0}^n x_j+(r-n)x_{n+1}$, where $n=\Int(r)$.
Thus, the above $d_i(s,m,\lambda)$'s are continuous functions in $s$, $m$ and
$\lambda$.

Note that
$\sum_{j=0}^{sm-1}(1+j-s\eta_\lambda(\frac js))(2+j-s\eta_\lambda(\frac js))/2=
\sum_{i=1}^3d_i(s,m,\lambda)s^i$ and \nl
$\lim_{s\rightarrow\infty}d_2(s,m,\lambda)=\int_0^m\frac32(x-\eta_\lambda(x))dx$,\qquad
$\lim_{s\rightarrow\infty}d_3(s,m,\lambda)=\int_0^m\frac12(x-\eta_\lambda(x))^2dx$.

(2.4) Now we divide the cases according to the values of $\lambda$.

1) The case $\lambda\le\frac73$.
We have $\eta_\lambda(\frac92)\le\frac92$.
Therefore
$h^0(M,sB)-h^0(M,s(\pi_1^*B-\frac92E))\le
\sum_{j=0}^{\frac92 s-1}\rk(\rho_{s,j})\le
\sum_{i=1}^3d_i(s,\frac92,\lambda)s^i$
by the above argument.
We have \nl
$\lim_{s\rightarrow\infty}6d_3(s,\frac 92,\lambda)=
{\dsize\int}_0^\frac92 3(x-\eta_\lambda(x))^2dx=
{\dsize\int}_0^{3-\frac 1\lambda}3x^2dx + {\dsize\int}_{3-\frac
1\lambda}^\frac92 3(x-\lambda(x-3)-1)^2dx=
\dfrac
1{\lambda-1}\left\{\dfrac{(3\lambda-1)^3}{\lambda^2}-\left(\dfrac{7-3\lambda}2\right)^3\right\}$.

Let $g(\lambda)$ be this function.
Then by computation we get
$g'(\lambda)=\frac{(3\lambda-1)}{4\lambda^3}h(\lambda)$
for $h(\lambda)=9\lambda^3-24\lambda^2-8\lambda+8$.
We have
$h'(\lambda)=27\lambda^2-48\lambda-8$ and $h''(\lambda)=54\lambda-48$.
In the interval $I=\{\frac13\le\lambda\le\frac73\}$, $h''$ changes the sign at
$\lambda=\frac89$, while $h'(\frac13)<0$ and $h'(\frac73)>0$.
Hence $h'(\lambda)$ changes the sign at some $\lambda_0$ with
$\frac89<\lambda_0<\frac73$.
Actually $\lambda_0=\frac{8+2\sqrt{22}}9$, but we do not need this explicit
value.

Any way we have $h(\frac13)=3>0$ and $h(\frac73)=-27<0$, so $h(\lambda_0)<0$.
Hence $h(\lambda)$ changes the sign at some $\lambda_1$ between $\frac13$ and
$\lambda_0$.
Therefore $g(\lambda)$ attains the maximum at $\lambda=\lambda_1$.

By computation we have $h(0.464)>0>h(0.465)$,
so $0.464<\lambda_1<0.465$, hence
$g(\lambda)<\dfrac1{1-0.465}\left\{\left(\dfrac{7-3\times0.464}2\right)^3-\dfrac{(3\times0.464-1)^3}{(0.465)^2}\right\}<40.69$.

Thus we have $\lim_{s\rightarrow\infty}6d_3(s,\frac92,\lambda)<41$ for any
$\frac13\le\lambda\le\frac73$, while $h^0(M,sB)$ grows like $ds^3/6$ with
$d=B^3\ge 51$.
Hence $h^0(M_1,s(\pi_1^*B-\frac92E))=O(s^3)$, i.e.,
$\kappa(\pi_1^*B-\frac92E)=3$ as desired.

2) The case $\lambda>\frac73$.
Set $T=\Inf\{t\in\BQ{\vert}\kappa(\pi_1^*B-tE)<0\}$
and $T_0=(3\lambda-1)/(\lambda-1)$.
Since $\xi(t)\le t$ if $\kappa(\pi_1^*B-tE)\ge 0$, we have
$T\le T_0<\frac92$,
since $\eta_\lambda(t)\le t$ iff $t\le T_0$.

We have $\lim_{s\rightarrow\infty}6d_3(s, T_0,\lambda)\le
\int_0^{3-\frac1\lambda}3x^2dx+\int_{3-\frac1\lambda}^{
T_0}3(x-\lambda(x-3)-1)^2dx$.
This right side $g(\lambda)$ is $(3\lambda-1)^3/(\lambda-1)\lambda^2$, so
$\frac d{d\lambda}\log g(\lambda)
=\frac9{3\lambda-1}-\frac1{\lambda-1}-\frac2\lambda
=\frac{-2}{(3\lambda-1)(\lambda-1)\lambda}<0$,
hence $g(\lambda)\le g(\frac73)<29.8$.
Therefore $\lim_{s\rightarrow\infty}6d_3(s,m,\lambda)<30$
for any $m$ such that $m- T_0$ is small enough.

If $ T\in\BQ$, set $m= T\le T_0$.
Then we have
$\sum_{j=0}^{sm-1}\rk(\rho_{s,j})\le\sum_{i=1}^3d_i(s,m,\lambda)s^i$ for any
$s\gg0$ and hence $\kappa(\pi_1^*B-mE)=3$.
But then $\kappa(\pi_1^*B-tE)>0$ for any $t$ slightly larger than $m= T$,
contradicting the choice of $ T$.

If $ T< T_0$, we take a rational number $m$ such that $ T<m< T_0$.
Then $\kappa(\pi_1^*B-mE)=3$ similarly as above, contradiction.

It remains the case $ T= T_0\notin\BQ$.
By similar computations as above, we infer that there are a small positive
number $\epsilon$ and a very large number $\ell$ such that
$h^0(M,sB)>\sum_{i=1}^3d_i(s,m,\lambda)s^i$ for any $m\le T+\epsilon$ and any
$s\le\ell$.
Here $d_i(s,m,\lambda)$ is continuous in $m$ under the convention in (2.3).
For an integer $s\ge\Max(\ell,\frac1\epsilon)$, let $q$ be the least integer
such that $q/s> T$, and set $m=q/s$.
Then
$\sum_{j=0}^{sm-1}\rk(\rho_{s,j})\le\sum_{i=1}^3d_i(s,m,\lambda)s^i<h^0(M,sB)$,
so $h^0(s\pi_1^*B-smE)>0$, contradicting $m> T$.

(2.5) Combining these observations, we conclude that $\pi_1^*B-\frac92E$ is big
unless we can get rid of non-exceptional places by replacing the degree of the
Hironaka decomposition as in (2.1).
\dnl
{\bf \S3. The case of exceptional place}

In this section we consider the cases where the locus of the bad place $P$ on
$E$ is a point.

(3.1) As a warm-up, we first consider the case where $P$ is the primary place
over a point $y$ on $E$.
Thus $P=\SP(E_1,Y_1)$ where $Y_1$ is the $(-1)$-curve on the blow-up $E_1$ of
$E$ at $y$.

(3.1.1) We have seen $\nu''_P-1\ge v_P(R_{F_0/E})=v_P(R_{E_1/E})=1$ in (1.7).
Hence, as in \S2, we can get rid of this situation by replacing the degree of
the Hironaka model when $\Xi_P(\pi_1^*B-3E)<2$.
Therefore we may assume that $\Xi_P(\pi_1^*B-3E)\ge2$.

(3.1.2) Set $\xi(t)=\Xi_P(\pi_1^*B-tE)$ for $t\in\BQ$.
Then, as in \S2, $\xi(t)$ is a convex function in $t$ with $\xi(3)\ge2$, so
$\xi(t)\ge\eta_\lambda(t)=_{\text{def}}\Max(\lambda(t-3)+2,0)$ for some
$\lambda\ge\frac23$.

(3.1.3) Set $h(u,w)=h^0(E_1,uH-wY_1)$, where $H$ is the pull-back of $\SO(1)$
on $E\cong\BP^2$.
Then we have $h(u,w)=\frac12(u+1)(u+2)-\frac12(w+1)w$ if $u\ge w$, and
$h(u,w)=0$ if $u<w$.

(3.1.4) Let $\rho_{s,j}: H^0(M_1,s\pi_1^*B-jE)\lra H^0(E,\SO_E(j))\cong
H^0(E_1,jH)$ be the natural map.
By the definition of $\Xi_P$, we have $\rk(\rho_{s,j})\le h(j,w)$ for some
$w\ge s\xi(\frac js)$, so \nl
$\rk(\rho_{s,j})\le\frac12(j+1)(j+2)-\frac12s\eta_\lambda(\frac
js)(s\eta_\lambda(\frac js)+1)$ if $\eta_\lambda(\frac js)\le \frac js$.\nl
Hence $h^0(M,sB)-h^0(M_1,s\pi_1^*B-smE)
\le\sum_{j=0}^{sm-1}\rk(\rho_{s,j})
\le\sum_{i=1}^3 d_i(s,m,\lambda)s^i$,
where we set \nl
$d_1(s,m,\lambda)=\sum_{j=0}^{sm-1}\frac 1s=m$,\nl
$d_2(s,m,\lambda)=\sum_{j=0}^{sm-1}(\frac32\frac js-\frac12\eta_\lambda(\frac
js))\frac1s$ and \nl
$d_3(s,m,\lambda)=\sum_{j=0}^{sm-1}(\frac12(\frac
js)^2-\frac12\eta_\lambda(\frac js)^2)\frac1s$.\nl
When $s\lra\infty$, we have
$d_2(s,m,\lambda)\lra\int_0^m(\frac32x-\frac12\eta_\lambda(x))dx$ and \nl
$d_3(s,m,\lambda)\lra\int_0^m(\frac12x^2-\frac12\eta_\lambda(x)^2)dx$.

(3.1.5) As in \S2, we divide the cases according to $\lambda$.

(3.1.5.1) Suppose that $\lambda\le\frac53$.
Then $\eta_\lambda(\frac92)\le\frac92$.
In this case we can apply the above argument for $m=\frac92$.
We have $\lim_{s\rightarrow\infty}6d_3(s,\frac92,\lambda)
=\int_0^\frac92 3(x^2-\eta_\lambda(x)^2)dx
=\int_0^{3-\frac2\lambda}3x^2dx +
\int_{3-\frac2\lambda}^{\frac92}3(x^2-(\lambda(x-3)+2)^2)dx
=(\frac92)^3-(\frac32\lambda+2)^3/\lambda$.
Set $g(\lambda)=(\frac32\lambda+2)^3/\lambda$.
Then $\frac d{d\lambda}\log g(\lambda)=\frac 9{3\lambda+4}-\frac1\lambda
=(6\lambda-4)/(3\lambda+4)\lambda>0$
in the range $\lambda>\frac23$.
Hence $g(\lambda)\ge g(\frac23)=\frac{81}2$ if $\lambda\ge\frac23$, so
$\lim_{s\rightarrow\infty}6d_3(s,\frac92,\lambda)\le(\frac92)^3-\frac{81}2=50.625<51\le d=B^3$.
By similar argument as in \S2, we infer $\kappa(\pi_1^*B-\frac92E)=3$ in this
case.

\Rmk The above number $50.625$ is the reason why we assume $d\ge51$ in the Main
Theorem.

(3.1.5.2) Suppose that $\lambda>\frac53$.
Then $\eta_\lambda( T_0)= T_0$ for $ T_0=(3\lambda-2)/(\lambda-1)<\frac92$.
We have $\int_0^{ T_0}3(x^2-\eta_{\lambda}(x)^2)dx
=(\frac{3\lambda-2}{\lambda-1})^3-\frac1\lambda(\frac{3\lambda-2}{\lambda-1})^3
=(3\lambda-2)^3/\lambda(\lambda-1)^2$, which we set $g(\lambda)$.
Then $\frac d{d\lambda}\log
g(\lambda)=\frac9{3\lambda-2}-\frac1\lambda-\frac2{\lambda-1}=-2/(3\lambda-2)\lambda(\lambda-1)<0$ if $\lambda>\frac53$,
so $g(\lambda)\le g(\frac53)=36.45<51$.
{}From this we get a contradiction as in (2.4.2).

(3.2) As a second warm-up, we consider the case where $P$ is a secondary place
over $y$.
This means that $P$ is the primary place over a point $y_2$ on the $(-1)$-curve
$Y_1$ on $E_1$.
Let $Y_2$ be the $(-1)$-curve over $y_2$ on the blow-up $E_2$ of $E_1$, so
$P=\SP(E_2,Y_2)$.

(3.2.1) In this case we have $\nu''_P-1\ge v_P(R_{E_2/E})=2$,
so we may assume $\Xi_P(\pi_1^*B-3E)\ge3$ as in (3.1.1).
Hence $\xi(t)=\Xi_P(\pi_1^*B-tE)$ satisfies
$\xi(t)\ge\eta_\lambda(t)=\Max(0,\lambda(t-3)+3)$ for some $\lambda\ge 1$.

(3.2.2) \Lem {\sl Let $H$ be the pull-back of $\SO(1)$ on $E\cong\BP^2$ to
$E_2$.
Set $h^i(u,w_1,w_2)=\dm H^i(E_2,uH-w_1Y_1^*-w_2Y_2)$, where $Y_1^*$ is the
total transform of $Y_1$.
Thus the proper transform $\tilde{Y_1}$ is a $(-2)$-curve on $E_2$ such that
$Y_1^*=\tilde{Y_1}+Y_2$.
Then\nl
1) $h^0(u,w_1,w_2)=\frac12(1+u)(2+u)-\sum_{i=1}^2\frac12w_i(w_i+1)$ if $u\ge
w_1+w_2$ and $w_1\ge w_2\ge 0$.\nl
2) $h^0(u,w_1,w_2)=0$ if $2u<w$, where $w=w_1+w_2$.\nl
3) $h^0(u,w_1,w_2)\le\frac12(1+u)(2+u)-\frac14w(w+2)$ if $w\le u$.\nl
4) $h^0(u,w_1,w_2)\le\frac14(2+2u-w)^2$ if $u\le w\le2u$.}

\Prf We claim $h^1(u,w_1,w_2)=0$ under the assumption in 1).
To prove this, we use the induction on $u$.
The claim is obvious if $w_2=0$, so suppose $u\ge w_1+w_2>0$.
The unique member $Y$ of ${\vert}H-Y_1^*-Y_2{\vert}$ is a $(-1)$-curve, so we
have an exact sequence \nl
$H^1(E_2,L(u-1,w_1-1,w_2-1))\lra H^1(E_2,L(u,w_1,w_2))\lra
H^1(Y,\SO_Y(u-w_1-w_2))$,\nl
where $L(u,w_1,w_2)$ denotes $uH-w_1Y_1^*-w_2Y_2\in\Pic(E_2)$.
The last term vanishes by the assumption $u\ge w_1+w_2$, while the first term
vanishes by the induction hypothesis.
Thus the claim is proved.

{}From this claim, 1) follows from the Riemann-Roch theorem.
2) is proved similarly by induction on $u$.
Indeed, we have an exact sequence \nl
$H^0(E_2,L(u-1,w_1-1,w_2-1))\lra H^0(E_2,L(u,w_1,w_2))\lra
H^0(Y,\SO(u-w_1-w_2))$,\nl
and the last term vanishes by the assumption $2u<w$.

To show 3), we may assume $w_1\ge w_2$.
Indeed, otherwise, we have an exact sequence \nl
$H^0(E_2,L(u,w_1+1,w_2-1))\lra
H^0(E_2,L(u,w_1,w_2))\lra
H^0(\tilde{Y_1},\SO(w_1-w_2))=0$, so
we reduce the problem to the case $w_1\ge w_2$ by the induction on $w_2-w_1$.
When $w_1\ge w_2$, we can apply 1), and we have
$\sum w_i(w_i+1)=\sum((w_i+\frac12)^2-\frac14)\ge
\frac12(\sum(w_i+\frac12))^2-\frac12=\frac12(w^2+2w)$ by Schwarz' inequality.
Combining them we get the estimate 3).

To prove 4), we use the induction on $w-u$.
We may assume $w_1\ge w_2$ as above.
The assertion follows from 3) for $w-u=0$.
When $w>u$, we use the exact sequence \nl
$H^0(E_2,L(u-1,w_1-1,w_2-1))\lra
H^0(E_2,L(u,w_1,w_2))\lra
H^0(Y,\SO(u-w_1-w_2))$.\nl
The last term vanishes by the assumption $w>u$, so the assertion is proved by
induction.

(3.2.3) Let $\rho_{s,j}$ be the map $H^0(M_1,s\pi_1^*B-jE)\lra
H^0(E,\SO(j))\cong H^0(E_2,jH)$.
By the definition of $\Xi$ we have $v_P(\psi)\ge s\xi(\frac js)$ for any
$\psi\in{\roman{Im}}(\rho_{s,j})$,
so  $\rk(\rho_{s,j})\le h^0(E_2,jH-wY_2)$ for some $w\ge s\xi(\frac js)$.
Hence, by (3.2.2), we have \nl
$\rk(\rho_{s,j})\le\frac12(1+j)(2+j)-\frac14s\eta_\lambda(\frac
js)(2+s\eta_\lambda(\frac js))$ if $\frac js\ge\xi(\frac js)$,\nl
\indent $\le\frac14(2+2j-s\eta_\lambda(\frac js))^2$ if $\frac js\le\xi(\frac
js)\le\frac{2j}s$ and \nl
\indent $=0$ if $\xi(\frac js)\ge\frac{2j}s$.\nl
According to these estimates, we modify the definition of $d_i(s,m,\lambda)$ in
(3.1.4) suitably, and argue as before.
This time we have
$$\lim_{s\rightarrow\infty}d_3(s,m,\lambda)=
\int_{\eta_\lambda(x)\le x}(\frac12
x^2-\frac14\eta_\lambda(x)^2)dx+\int_{x\le\eta_\lambda(x)\le2x}\frac14(2x-\eta_\lambda(x))^2dx$$.

(3.2.4) As in (3.1.5), we divide the cases according to $\lambda$.

(3.2.4.1) $\lambda\le4$.
In this case the above limit is
$\int_0^{3-\frac3\lambda}\frac12 x^2dx +
\int_{3-\frac3\lambda}^3\{\frac12x^2-\frac14(\lambda(x-3)+3)^2\}dx\almb
+\int_3^\frac92\frac14(2x-\lambda(x-3)-3)^2dx
=\frac92-\frac9{4\lambda}+\frac9{32}(\lambda^2-10\lambda+28)
=\frac9{32}(\lambda^2-10\lambda+44-\frac8\lambda)$
for $m=\frac92$.
It is easy to see that this is a decreasing function in $\lambda$, so
$\le\frac{243}{32}=7.59375<\frac{51}6$ in the range $1\le\lambda\le4$.

(3.2.4.2) $\lambda\ge4$.
In this case $\eta_\lambda(x)\le2x$ only if $x\le
T_0=(3\lambda-3)/(\lambda-2)$, and the above limit for $m= T_0$ is
$\frac92(1+\frac1{(\lambda-2)\lambda})$.
This is a decreasing function in $\lambda$ when $\lambda\ge4$, so
${}\le\frac{81}{16}=5.0625$.

(3.2.5) In either case  we have the limit ${}<\frac{51}6$, so we are done as in
\S2 or (3.1).

(3.3) Now we consider the general case.
By succesive blow-ups at the loci of $P$, we get a model $\tilde E$ of $E$
together with a divisor $D$ on it such that $P=\SP(\tilde E,D)$.
$(\tilde E,D)$ will be called the lowest representative pair of $P$ and we will
divide the cases according to its nature.

The blow-up process is essentially the same up to the second step, where the
situation is as in (3.2).
The locus of $P$ on $E_2$, the center of the third blow-up, is a point $y_3$ on
the newest $(-1)$-curve $Y_2$.
There are now two cases: $y_3=Y_2\cap\tilde{Y_1}$ or $y_3\notin\tilde{Y_1}$.
In the former case, the $(-1)$-curve $Y_3$ over $y_3$ meets the proper
transforms of the other older $(-1)$-curves at two points, while in the latter
case $Y_3$ meets only $\tilde Y_2$ at a point.
Thus, also in the subsequent steps, the locus of $P$ is a point on the newest
$(-1)$-curve, which meets the proper transforms of older $(-1)$-curves at one
or two point(s).

The whole process can be described as follows.
At first there are several times, say $i$, of blow-ups at points off the proper
transforms of older $(-1)$-curves.
They form a chain as below, where $(n)$ denotes the proper transform of the
$n$-th $(-1)$-curve.
$$(1)-(2)-\cdots-(i)\tag1$$
Next come blow-ups at the meeting point with the proper transform of the oldest
$(-1)$-curve.
Call the proper transforms of the $(-1)$-curves of these steps $(i.2)$,
$(i.3)$, $\cdots$, $(i.j)$ and rename $(i)$ as $(i.1)$.
The result is a chain as below:
$$(1)-\cdots-(i-1)-(i.j)-(i.j-1)-\cdots-(i.1)\tag2$$
Next come blow-ups at the meeting points with $(i.j-1)$ above.
Renaming $(i.j)$ as $(i.j.1)$, we get a chain as follows:
$$(1)-\cdots-(i-1)-(i.j.1)-\cdots-(i.j.k)-(i.j-1)-\cdots-(i.1)\tag3$$
Similar process will continue unless we blow up at a point off the proper
transforms of older $(-1)$-curves.
For example, the next step of the above will be
$$(1)-\cdots-(i-1)-(i.j.1)-\cdots-(i.j.k-1)-(i.j.k.\ell)-\cdots-(i.j.k.1)-(i.j-1)-\cdots-(i.1)\tag4$$
Note that $i$, $j$, $k$, $\ell$, $\cdots$ are integers $\ge2$.

If we blow-up at a point off the proper transforms of older $(-1)$-curves, at
first the situation is as follows:
$$\vdash-[2]\tag5$$
Here $\vdash$ is the union of proper transforms of the older $(-1)$-curves, and
[2] is the newest $(-1)$-curve.
Then comes a similar process as before, and a sample result will be:
$$\vdash-[2]-\cdots-[i'-1]-[i'.j'.1]-\cdots-[i'.j'.k']-[i'.j'-1]-\cdots-[i'.1]\tag6$$
Now comes a blow-up at a point off the proper transforms of older
$(-1)$-curves, followed by a similar process as above, and so on.

After several times we get a pair representing $P$.

(3.4) Let $\Lambda$ be the restriction to $E$ of the linear system ${\vert}\tau
(\pi_1^*B-3E){\vert}$, where $\tau $ is the degree of the Hironaka model.
Set $m_1=\tau^{-1}v_{Y_1}(\Lambda)$.
Then $\Lambda_1=\Lambda{\vert}_{E_1}-\tau m_1Y_1$ on $E_1$ has no exceptional
fixed component.
Next set $m_2=\tau ^{-1}v_{Y_2}(\Lambda_1)$.
Then $\Lambda_2=\Lambda_1{\vert}_{E_2}-\tau m_2Y_2$ on $E_2$ has no exceptional
fixed component.
We proceed similarly and define $m_*$ and $\Lambda_*$.
For example, if $E_{(i.j.k)}$ is the intermediate step as in (3.3;3) and if
$Y_{(i.j.k)}$ is the newest $(-1)$-curve on it, set
$m_{(i.j.k)}=\tau ^{-1}v_{Y_{(i.j.k)}}(\Lambda_{(i.j.k-1)})$.
Then $\Lambda_{(i.j.k)}=\Lambda_{(i.j.k-1)}{\vert}_{E_{(i.j.k)}}-\tau
m_{(i.j.k)}Y_{(i.j.k)}$ on $E_{(i.j.k)}$ has no exceptional fixed component.
In the same way $m_{[i'.j']}$ etc. are defined.

(3.5) The above numbers $\{m_*\}$ satisfy various inequalities.
Suppose for example we are in the situation (3) in (3.3):
$$(1)-\cdots-(i-1)-(i.j.1)-\cdots-(i.j.k)-(i.j-1)-\cdots-(i.1)$$
Then we have\nl
(1) $3\ge m_{(1)}\ge\cdots\ge m_{(i-1)}\ge
m_{(i.1)}+\cdots+m_{(i.j-1)}+m_{(i.j.1)}$.

Indeed, for the proper transform $\tilde Y_{(*)}$ of $Y_{(*)}$, we have
$0\le\Lambda_{(*+1)}\{\tilde Y_{(*)}\}=m_{(*)}-m_{(*+1)}$ for each $*=1,
\cdots, i-2$.
For $*=i-1$ we get the last inequality from
$0\le\Lambda_{(i.j.1)}\{Y_{(i-1)}\}$.

Similarly we obtain\nl
(2) $m_{(i.1)}\ge\cdots\ge m_{(i.j-1)}\ge m_{(i.j.1)}+\cdots+m_{(i.j.k)}$.\nl
(3) $m_{(i.j.1)}\ge\cdots\ge m_{(i.j.k)}$.

In the situation (3.3;5), we have $m_0\ge m_{[2]}$, where $m_0$ is the $m$ of
the $(-1)$-curve just before $Y_{[2]}$.
In case (3.3;6) we have similarly \nl
$m_0\ge m_{[2]}\ge\cdots\ge m_{[i'-1]}\ge
m_{[i'.1]}+\cdots+m_{[i'.j'-1]}+m_{[i'.j'.1]}$.

Probably these samples are enough to understand the general principle.

(3.6) Set $r_*=v_*(R_*)$ for $*=(i), (i.j), \cdots, [i'], \cdots$ etc.,
where $R_*$ is the ramification divisor of the map $E_*\lra E$.
Then $r_*-1$ is the sum of $r_{\#}$'s such that the blow-up center $y_*$ lies
on the proper transform of $Y_{\#}$.
In particular, we have \nl
$r_{(n)}=n$ for $n\le i$,\nl
$r_{(i.n)}=i+(n-1)(r_{(i-1)}+1)=ni$ for $n\le j$,\nl
$r_{(i.j.n)}=r_{(i.j)}+(n-1)(r_{(i.j-1)}+1)=nij-(n-1)(i-1)$ for $n\le k$,\nl
$r_{(i.j.k.n)}=r_{(i.j.k)}+(n-1)(r_{(i.j.k-1)}+1)=ijkn-ijn-ikn+ij+nk+2ni-n-i$,\nl
$\cdots\cdots$\nl
and similarly we can calculate $r_*$ in the general case.

Next set $\mu_*=v_*(\Lambda{\vert}_{E_*})$.
Then $\mu_*-m_*$ is the sum of $\mu_{\#}$'s such that the blow-up center $y_*$
lies on the proper transform of $Y_{\#}$.
In particular, we have\nl
$\mu_{(n)}=m_{(1)}+\cdots+m_{(n)}$ for $n\le i$,\nl
$\mu_{(i.n)}=n\mu_{(i-1)}+m_{(i.1)}+m_{(i.2)}+\cdots+m_{(i.n)}$\nl
\indent $=n(m_{(1)}+\cdots+m_{(i-1)}+m_{(i.1)}+\cdots+m_{(i.n)}$ for $n\le
j$,\nl
$\mu_{(i.j.n)}=n\mu_{(i.j-1)}+\mu_{(i-1)}+m_{(i.j.1)}+\cdots+m_{(i.j.n)}$\nl
\indent
$=((j-1)n+1)(m_{(1)}+\cdots+m_{(i-1)})+n(m_{(i.1)}+\cdots+m_{(i.j-1)})+m_{(i.j.1)}+\cdots+m_{(i.j.n)}$ for $n\le k$,\nl
$\mu_{(i.j.k.n)}=\mu_{(i.j-1)}+n\mu_{(i.j.k-1)}+m_{(i.j.k.1)}+\cdots+m_{(i.j.k.n)}$\nl
\indent
$=((j-1)+n((j-1)(k-1)+1)(m_{(1)}+\cdots+m_{(i-1)})+(n(k-1)+1)(m_{(i.1)}+\cdots+m_{(i.j-1)})+n(m_{(i.j.1)}+\cdots+m_{(i.j.k-1)})+m_{(i.j.k.1)}+\cdots+m_{(i.j.k.n)}$,\nl
$\cdots\cdots$

Now we set $\delta_*=\mu_*-r_*$, which can be calculated as above.
Note that $\mu_*$ corresponds $\nu''_*$ in (1.7), hence $\delta_*\ge 1$ iff
$P_*$ is a bad place.
Hence we may assume $\delta_*<1$ for any $*$ except the final one.
In particular we may assume $m_{(1)}<2$.

(3.7) \Lem {\sl Set
$\Delta_{(i.1)}=\delta_{(i.1)}+m_{(i-1)}$,
$\Delta_{(i.j.1)}=\delta_{(i.j.1)}+m_{(i.j-1)}$,
$\Delta_{(i.j.k.1)}=\delta_{(i.j.k.1)}+m_{(i.j.k-1)}$,
$\cdots$,
$\Delta_{[i'.1]}=\delta_{[i'.1]}+m_{[i'-1]}$,
$\cdots$, etc.
Then
$\Delta_{(i.1)}\ge\Delta_{(i.j.1)}\ge\Delta_{(i.j.k.1)}\ge\cdots\ge\Delta_{[i'.1]}\ge\Delta_{[i'.j'.1]}\ge\cdots\ge1$.}

\Prf We have
$\delta_{(i.j.1)}=\delta_{(i.1)}+(j-1)(\delta_{(i-1)}-1)+m_{(i.2)}+\cdots+m_{(i.j)}
\le \delta_{(i.1)}+m_{(i.2)}+\cdots+m_{(i.j)}$.
Hence $\Delta_{(i.j.1)}\le\delta_{(i.j.1)}+m_{(i.1)}\le\Delta_{(i.1)}$ by
(3.5).
Similarly we have $\delta_{(i.j.k.1)}\le
\delta_{(i.j.1)}+m_{(i.j.2)}+\cdots+m_{(i.j.k)}$ and
$\Delta_{(i.j.k.1)}\le\delta_{(i.j.k.1)}+m_{(i.j.1)}\le\Delta_{(i.j.1)}$, and
so on.

In the situation (6) in (3.3), let $Y_0$ be the newest $(-1)$-curve in the part
$\vdash$.
Then $\Delta_0\ge\delta_0+m_0$, $m_0\ge m_{[i'-1]}$ and
$\delta_{[i'.1]}=\delta_0+\sum_n(m_{[n]}-1)\le\delta_0$,
so $\Delta_0\ge\Delta_{[i'.1]}$.

By similar argument we obtain the asserted inequalities.
As for the final one, let us consider the case where $Y_{[i'.j'.k']}$ is the
final $(-1)$-curve as in (3.3.6), as a sample.
We have
$1\le\delta_{[i'.j'.k']}=\delta_{[i'.j'.1]}+(k'-1)(\delta_{[i'.j'-1]}-1)+m_{[i'.j'.2]}+\cdots+m_{[i'.j'.k']}\le\delta_{[i'.j'.1]}+\sum_n m_{[i'.j'.n]}\le\Delta_{[i'.j'.1]}$.
Other cases can be treated similarly.

\Rmk If you define $\Delta_{(i.n)}=\delta_{(i.n)}+m_{(i.n-1)}$ for $n>1$, then
this is not always a decreasing function in $n$.

(3.8) Let $\Xi_*(\pi_1^*B-3E)$ (sometimes abbreviated as $\Xi_*$ from now on)
be the algebraic Lelong number at the place $P_*=\SP(Y_*)$.
By definition, this is approximated by $\mu_*$ when the degree $\tau $ of the
Hironaka model is large enough.
In the sequel we will show that the bigness of $\pi_1^*B-\frac92E$ follows from
lower estimates of $\Xi_*$ as in (3.1) and (3.2).
But the precise arguments depend on the place $P_*$.

(3.9) First we consider the case $P_{(n)}$, $1\le n\le i$.
Let $H$ be the pull-back of $\SO(1)$ on $E\cong\BP^2$ to $E_{(n)}$ and let
$Y_j$ be the total transform of the $(-1)$-curve $Y_{(j)}$ on $E_{(n)}$, while
the proper transform is denoted by $\tilde{Y_j}$.
Thus, $\tilde{Y_j}$ is the unique member of ${\vert}Y_j-Y_{j+1}{\vert}$ for
each $j<n$.
Put $L(u, w_1, \cdots, w_n)=uH-w_1Y_1-\cdots-Y_n\in\Pic(E_{(n)})$,
$h^p(u, w_1, \cdots, w_n)=\dm H^p(E_{(n)}, L(u,w_1,\cdots,w_n))$ and
$w=w_1+\cdots+w_n$ for non-negative integers $u, w_1, \cdots, w_n$.

Then we have the following

\Lem {\sl
1) $h^0(u,w_1,\cdots,w_n)=0$ if $nu<w$.\nl
2) $h^1(u,w_1,\cdots,w_n)=0$ and
$h^0(u,w_1,\cdots,w_n)=\frac12(u+1)(u+2)-\frac12\sum_j w_j(w_j+1)$ if $u\ge w$
and if $w_1\ge w_2\ge \cdots \ge w_n$.\nl
3) $h^0(u,w_1,\cdots,w_n)\le\frac12(u+1)(u+2)-\frac1{2n}w(w+n)$ if $u\ge w$.\nl
4) $h^0(u,w_1,\cdots,w_n)\le q(nu-w)$ if $u\le w\le nu$, where
$q(\theta)=\dfrac{\theta^2}{2n(n-1)}+\dfrac\theta{n-1}+1$.
}

\Prf. The point $y_{(2)}$ on $Y_1$ determines a line $Z$ on $E$ passing the
point $y_{(1)}$.
Its proper transform $\tilde Z$ on $E_{(n)}$ is a member of
${\vert}H-Y_1-\cdots-Y_a{\vert}$ for some $a$ with $2\le a\le n$.

To show 1), we use the induction on $u$.
The assertion is obvious when $u=0$, so suppose $u\ge1$.
If $w_j<w_{j+1}$ for some $j$, then $L(u,w_1,\cdots,w_n)\tilde{Y_j}<0$, so
$h^0(u,w_1,\cdots,w_n)\le h^0(u,w_1,\cdots,w_j+1,w_{j+1}-1,\cdots)$
since \nl
$H^0(L(u,w_1,\cdots,w_n)-\tilde{Y_j})\lra H^0(L(u,w_1,\cdots,w_n))\lra
H^0(\tilde{Y_j},L(u,w_1,\cdots,w_n){\vert}_{Y_j})$ \nl
is exact.
Therefore it suffices to consider the case $w_1\ge \cdots\ge w_n$.
Then $L(u,w_1,\cdots,w_n)\tilde Z=u-w_1-\cdots-w_a\le u-a\frac wn<0$, so we
have
$h^0(u,w_1,\cdots,w_n)=h^0(L(u,w_1,\cdots,w_n)-\tilde
Z=h^0(u-1,w_1-1,\cdots,w_a-1,w_{a+1},\cdots,w_n)$.
This last term is zero by the induction hypothesis, thus 1) is proved.

We use the induction on $u$ to prove 2) too.
It suffices to prove $h^1=0$, since the assertion on $h^0$ follows from this by
the Riemann-Roch theorem.
The case $u=0$ is obvious, so suppose $u\ge 1$.
By assumption $L(u,w_1,\cdots,w_n)\tilde Z=u-w_1-\cdots-w_a\ge0$,
so $h^1(u,w_1,\cdots,w_n)\le h^1(L(u,w_1,\cdots,w_n)-\tilde
Z)=h^1(u-1,w_1-1,\cdots,w_a-1,w_{a+1},\cdots)$.
If $w_a>w_{a+1}$, this vanishes by the induction hypothesis.
If $w_a=w_{a+1}$, we have $L(u-1,\cdots,w_a-1,w_{a+1},\cdots)\tilde{Y_a}=-1$
and $H^1(\tilde{Y_a},\SO(-1))=0$, hence
$h^1(u-1,w_1-1,\cdots,w_a-1,w_{a+1},\cdots)\le h^1(u-1,\cdots,w_{a-1}-1,\almb
w_a,w_{a+1}-1,w_{a+2},\cdots)$.
Now if $w_{a-1}-1<w_a$ ($\Leftrightarrow w_{a-1}=w_a$), we proceed similarly by
subtracting $\tilde Y_{a-1}$ to obtain
$\le h^1(u-1,\cdots,w_{a-2}-1,w_{a-1},w_a,w_{a+1}-1,w_{a+2},\cdots)$.
After several similar steps we get
$\le h^1(u-1,w_1 (\text{or } w_1-1),\cdots,w_{a+1}-1,w_{a+2},\cdots)$
where we have $w_1 (\text{or } w_1-1)\ge\cdots\ge w_{a+1}-1$.
If $w_{a+1}>w_{a+2}$ we are done by the induction hypothesis.
If $w_{a+1}=w_{a+2}$, subtracting $\tilde{Y_j}$'s several times as above,
we get $\le h^1(u-1,w_1 (\text{or }
w_1-1),\almb\cdots,w_{a+2}-1,w_{a+3},\cdots)$
such that $w_1 (\text{or } w_1-1)\ge\cdots\ge w_{a+2}-1$.
Repeating this process we get $\le h^1(u-1,w'_1,\cdots,w'_n)$ such that
$w'_1\ge\cdots\ge w'_n$ and $w'_1+\cdots+w'_n=w-a\ge w-n$.
This $h^1=0$ by the induction hypothesis, thus 2) is proved.

3) follows from 2).
Indeed, $\sum w_j(w_j+1)=\sum(w_j+\frac12)^2-\frac
n4\le\frac1n(\sum(w_j+\frac12))^2-\frac n4=\frac1n(w+\frac n2)^2-\frac
n4=\frac1n w^2+w$
by Schwarz inequality.

Note that both 2) and 3) are valid under the weaker assumption $u\ge
w_1+\cdots+w_a$ instead of $u\ge w$.

We use the induction on $u$ to prove 4).
The assertion is trivial for $u=0$, so suppose $u>0$.
If $w_j<w_{j+1}$ for some $j$, we have
$h^0(u,w_1,\cdots,w_n)=h^0(u,\cdots,w_j+1,w_j-1,\cdots)$
since $H^0(\tilde{Y_j},L(u,w_1,\cdots,w_n))=0$.
Therefore, as in the proof of 2), we reduce the problem to the case
$w_1\ge\cdots\ge w_n$.
Now if $u\ge w_1+\cdots+w_a$, then 3) applies.
Since
$\frac12(u+1)(u+2)\almb-\frac1{2n}w(w+n)-q(nu-w)=\frac1{2(n-1)}(n-3-u+w)(u-w)\le 0$, we are done in this case.
If $u<w_1+\cdots+w_a$, we have $L(u,\cdots){\tilde Z}<0$.
Therefore
$h^0(u,\cdots)=h^0(L(u,\cdots)-\tilde Z)=h^0(u-1,\almb
w_1-1,\cdots,w_a-1,w_{a+1},\cdots)\le q(nu-w-n+a)$ by the induction hypothesis.
Since $q(\theta-n+a)-q(\theta)=-\frac{n-a}{2n(n-1)}(2\theta+n+a)<0$,
we have $h^0(u,\cdots)\le q(nu-w)$ in this case too, as desired.

(3.10) \Lem {\sl $\pi_1^*B-\frac92 E$ is big in any of the following cases:\nl
1) $\Xi_{(1)}>\sqrt{321}/9=1.990719\cdots$\nl
2) $\Xi_{(2)}>\sqrt{642}/9=2.81530\cdots$\nl
3) $\Xi_{(3)}\ge3.51$ .\nl
4) $\Xi_{(4)}\ge 4.23$ .\nl
5) $\Xi_{(5)}\ge 5$ .}

For the proof, we use the following computational result, which is more general
than is needed here, and is used later repeatedly.

(3.11) \Lem {\sl For constants $\alpha$, $\beta$ and $\gamma$, put
$$f_{\alpha,\beta,\gamma}(\lambda)=\frac{(3\lambda-\gamma)^3}{\lambda(\lambda-\alpha)(\lambda-\beta)}-\frac{(9\beta-2\gamma-3\lambda)^3}{8\beta(\beta-\alpha)(\lambda-\beta)} {\text {\qquad and\qquad}}\epsilon=\alpha+\beta-\gamma .$$
Then
$$\frac
d{d\lambda}f_{\alpha,\beta,\gamma}(\lambda)=\frac{3\lambda-\gamma}{4\beta(\beta-\alpha)\lambda^2(\lambda-\alpha)^2}\varphi(\lambda)$$
where
$\varphi(\lambda)=9\lambda^4-(18\alpha+36\beta-12\gamma)\lambda^3+(9\alpha(\alpha+8\beta)-12(2\alpha+\beta)\gamma+4\gamma^2)\lambda^2\almb+4(\alpha-\beta)\gamma(3\alpha-2\gamma)\lambda\almb+4\alpha(\alpha-\beta)\gamma^2\almb
=9\lambda^4-6(\alpha+4\beta+2\epsilon)\lambda^3\almb+((-11\alpha^2+44\alpha\beta-8\beta^2)+(16\alpha+4\beta)\epsilon+4\epsilon^2)\lambda^2\almb+4(\beta-\alpha)((-\alpha^2+\alpha\beta+2\beta^2)-(\alpha+4\beta)\epsilon+2\epsilon^2)\lambda\almb+4\alpha(\alpha-\beta)((\alpha+\beta)^2-2(\alpha+\beta)\epsilon+\epsilon^2)$.\qquad
Moreover
$$\align
\varphi(\tfrac\gamma3) &=\gamma^2(\gamma-3\alpha)^2,\\
\varphi(3\alpha-\tfrac23\gamma)
&=36\alpha(\alpha-\beta)(2\alpha-\beta+\epsilon)^2,\\
\varphi(3\beta-\tfrac23\gamma)
&=9\beta(\alpha-\beta)((4\alpha^2-7\alpha\beta+7\beta^2)-8(\alpha-2\beta)\epsilon+4\epsilon^2),\\
\varphi'(\tfrac\gamma3) &=2\gamma(\gamma-3\alpha)(4\gamma-3\alpha-6\beta),\\
\varphi'(3\alpha-\tfrac23\gamma)
&=12(\alpha-\beta)(13\alpha-2\beta+2\epsilon)(2\alpha-\beta+\epsilon),\\
\varphi'(3\beta-\tfrac23\gamma)
&=18(\alpha-\beta)\beta(\alpha-2\beta+2\epsilon),\\
\varphi''(\tfrac\gamma3)
&=44\gamma^2-12(7\alpha+8\beta)\gamma+18\alpha(\alpha+8\beta) \\
&=(-22\alpha^2+52\alpha\beta-52\beta^2)+(-4\alpha+8\beta)\epsilon+44\epsilon^2,\\
\varphi''(3\alpha-\tfrac23\gamma)
&=482\alpha^2-560\alpha\beta+128\beta^2+(176\alpha-136\beta)\epsilon+8\epsilon^2,\\
\varphi''(3\beta-\tfrac23\gamma)
&=50\alpha^2-236\alpha\beta+236\beta^2-40(\alpha-2\beta)\epsilon+8\epsilon^2,\\
\varphi^{'''}(\tfrac\gamma3) &=36(-3\alpha-6\beta+4\gamma),\\
\varphi^{'''}(3\alpha-\tfrac23\gamma) &=36(13\alpha-8\beta+2\epsilon),\\
\varphi^{'''}(3\beta-\tfrac23\gamma) &=36(-5\alpha+10\beta+2\epsilon).
\endalign$$}

\Prf By elementary computation.

(3.12) Proof of (3.10).
The cases 1) and 2) are treated by the same method as in (3.1) and (3.2).
For example, in case 2), let $\eta_\lambda(x)=\Max(0,\lambda(x-3)+\gamma)$ with
$\gamma>\sqrt{642}/9$.
In place of (3.2.4.1), we should show
$$\int_0^{3-\frac{\gamma}\lambda}3x^2dx+\int_{3-\frac
\gamma\lambda}^{\frac92}(3x^2-\frac32\eta_\lambda(x)^2)dx=\frac{729}8-\frac1{16\lambda}(3\lambda+2\gamma)^3<51$$
for $\frac \gamma3\le\lambda\le3-\frac23\gamma$ and
$$\multline
\int_0^{3-\frac \gamma\lambda}3x^2dx+\int_{3-\frac
\gamma\lambda}^{\frac{3\lambda-\gamma}{\lambda-1}}(3x^2-\frac32\eta_\lambda(x)^2)dx+\int_{\frac{3\lambda-\gamma}{\lambda-1}}^\frac92\frac32(2x-\eta_\lambda(x))^2dx\\
=\frac{(3\lambda-\gamma)^3}{\lambda(\lambda-1)(\lambda-2)}-\frac{(18-3\lambda-2\gamma)^3}{16(\lambda-2)}<51
\endmultline$$
for $3-\frac23\gamma\le\lambda\le6-\frac23\gamma$.
It is easy to see that $(3\lambda+2\gamma)^3/\lambda$ is an increasing function
for $\lambda>\frac \gamma3$, so the assertion is valid for $\lambda\le3-\frac23
\gamma$.
In case $\lambda\ge3-\frac23 \gamma$, we apply (3.11) for $\alpha=1$, $\beta=2$
and $\epsilon=3-\gamma=0.1847\cdots$.
We have $\varphi^{'''}(3-\frac23\gamma)=36(-3+2\epsilon)<0$ and
$\varphi^{'''}(6-\frac23\gamma)=36(15+2\epsilon)>0$, so $\varphi''$ decreases
first and then increases in the interval $I=\{3-\frac23\gamma\le
\lambda\le6-\frac23\gamma\}$.
Since $\varphi''(3-\frac23\gamma)=-126-96\epsilon+8\epsilon^2<0$ and
$\varphi''(6-\frac23\gamma)=522+120\epsilon+8\epsilon^2>0$,
the sign of $\varphi''$ changes only once in $I$ from $-$ to $+$.
We have $\varphi'(3-\frac23\gamma)=-12(9+2\epsilon)\epsilon<0$ and
$\varphi'(6-\frac23\gamma)=36(3-2\epsilon)>0$,
so the sign of $\varphi'$ changes once from $-$ to $+$.
We have $\varphi(3-\frac23\gamma)=-36\epsilon^2<0$ and
$\varphi(6-\frac23\gamma)=-18(18+24\epsilon+4\epsilon^2)<0$,
so $\varphi(\lambda)<0$ for any $\lambda\in I$, hence
$f_{\alpha,\beta,\gamma}(\lambda)$ is decreasing.
This yields the desired assertion.

In place of (3.2.4.2), we should show
$$\int_0^{3-\frac \gamma\lambda}3x^2dx+\int_{3-\frac
\gamma\lambda}^{\frac{3\lambda-\gamma}{\lambda-1}}(3x^2-\frac32\eta_\lambda(x)^2)dx+\int_{\frac{3\lambda-\gamma}{\lambda-1}}^{\frac{3\lambda-\gamma}{\lambda-2}}\frac32(2x-\eta_\lambda(x))^2dx=\frac{(3\lambda-\gamma)^3}{\lambda(\lambda-1)(\lambda-2)}<51$$
for $\lambda\ge6-\frac23\gamma$.
This is easy since
$(3\lambda-\gamma)\lambda(\lambda-1)(\lambda-2)\dfrac
d{d\lambda}\log(\dfrac{(3\lambda-\gamma)^3}{\lambda(\lambda-1)(\lambda-2)})\almb
=((3\gamma-9)\lambda^2+(12-6\gamma)\lambda+2\gamma)
=-\epsilon\lambda^2+(-6+6\epsilon)\lambda+2(3-\epsilon)=6(1-\lambda)-\epsilon(\lambda^2-6\lambda+2)<0$.

The other cases 3), 4) and 5) can be treated similarly.
In case 5), we have $\xi_{(5)}(x)\ge\eta_\lambda(x)=\Max(0,\lambda(x-3)+5)$ for
some $\lambda\ge\frac53$.
By applying (3.9) as in (3.2.3), we get \nl
$\rk(\rho_{s,j})\le\frac12(1+j)(2+j)-\frac1{10}s\eta_\lambda(\frac
js)(s\eta_\lambda(\frac js)+5)$ if $\frac js\ge\xi_{(5)}(\frac js)$,\nl
$\rk(\rho_{s,j})\le\frac1{40}(5j-s\eta_\lambda(\frac
js))^2+\frac14(5j-s\eta_\lambda(\frac js))+1$ if $\frac js\le\eta_\lambda(\frac
js)\le\frac{5j}s$ and \nl
$\rk(\rho_{s,j})=0$ if $\eta_\lambda(\frac js)>\frac{5j}s$.\nl
Hence we should estimate
$$\lim_{s\rightarrow\infty}d_3(s,m,\lambda)=\int_{\eta_\lambda(x)\le
x}(\frac12x^2-\frac1{10}\eta_\lambda(x)^2)dx+\int_{x\le\eta_\lambda(x)\le5x}\frac1{40}(5x-\eta_\lambda(x))^2dx.$$
As before, for $\gamma=5$, we show
$$\multline
\int_0^{3-\frac \gamma\lambda}3x^2dx+\int_{3-\frac
\gamma\lambda}^{\frac{3\lambda-\gamma}{\lambda-1}}(3x^2-\frac35\eta_\lambda(x)^2)dx+\int_{\frac{3\lambda-\gamma}{\lambda-1}}^\frac92\frac6{40}(5x-\eta_\lambda(x))^2dx\\
=\frac{(3\lambda-\gamma)^3}{\lambda(\lambda-1)(\lambda-5)}-\frac{(45-2\gamma-3\lambda)^3}{160(\lambda-5)}<51
\endmultline$$
for $\frac \gamma3\le\lambda\le15-\frac23\gamma$ and
$$\int_0^{3-\frac
\gamma\lambda}3x^2dx+\cdots+\int_{\frac{3\lambda-\gamma}{\lambda-1}}^{\frac{3\lambda-\gamma}{\lambda-5}}\frac6{40}(5x-\eta_\lambda(x))^2dx=(3\lambda-\gamma)^3/\lambda(\lambda-1)(\lambda-5)<51$$
for $\lambda\ge15-\frac23\gamma$.
This latter case is easy.
In the former case we apply (3.11) for $\alpha=1$, $\beta=5$, $\gamma=5$ and
$\epsilon=1$.
In the interval $I=\{\frac \gamma3\le
\lambda\le15-\frac23\gamma\}=[\frac53,\frac{35}3]$, $\varphi^{'''}$ increases
from $-$ to $+$.
Since $\varphi''(\frac53)<0$ and $\varphi''(\frac{35}3)>0$, the sign of
$\varphi''$ changes once from $-$ to $+$.
We have $\varphi'(\frac53)=-260<0$ and $\varphi'(\frac{35}3)=2520>0$, so
$\varphi'$ changes its sign from $-$ to $+$.
Since $\varphi(\frac53)=100>0$ and $\varphi(\frac{35}3)=-39600<0$, the sign of
$\varphi$ changes from $+$ to $-$ at some $\lambda_1\in I$, and
$f_{\alpha,\beta,\gamma}(\lambda)$ attains the maximum at $\lambda=\lambda_1$.
By computation we see $\varphi(1.922)>0$ and $\varphi(1.923)<0$, hence
$1.922<\lambda_1<1.923$, so
$f(\lambda_1)<\dfrac1{5-1.923}\left(\dfrac{(35-3\times1.922)^3}{160}-\dfrac{(3\times1.922-5)^3}{1.923\times0.923}\right)=50.66\cdots<51$, as desired.

In case 4), we have
$\xi_{(4)}(x)\ge\eta_\lambda(x)=\Max(0,\lambda(x-3)+\gamma)$ for $\gamma=4.23$
for some $\lambda\ge\frac \gamma3$.
Similarly as before, we have
$$\lim_{s\rightarrow\infty}d_3(s.m.\lambda)=\int_{\eta_\lambda(x)\le
x}(\frac12x^2-\frac18\eta_\lambda(x)^2)dx+\int_{x\le\eta_\lambda(x)\le4x}\frac1{24}(4x-\eta_\lambda(x))^2dx .$$
We show
$\dfrac{(3\lambda-\gamma)^3}{\lambda(\lambda-1)(\lambda-4)}-\dfrac{(36-2\gamma-3\lambda)^3}{96(\lambda-4)}<51$
for $\frac \gamma3<\lambda\le12-\frac23\gamma$ and
$\dfrac{(3\lambda-\gamma)^3}{\lambda(\lambda-1)(\lambda-4)}<51$ for
$\lambda\ge12-\frac23\gamma$.
The latter is easy.
In the former case we apply (3.11) for $\alpha=1$, $\beta=4$, $\gamma=4.23$ and
$\epsilon=0.77$.
Then
$\varphi(\lambda)=9\lambda^4-111.24\lambda^3+64.0116\lambda^2+277.1496\lambda-214.7148$ and
$\varphi(1.581)>0>\varphi(1.582)$, so
$f(\lambda)$ attains the maximum at $\lambda=\lambda_1$ with
$1.581<\lambda_1<1.582$.
Hence
$f(\lambda_1)<\dfrac1{4-1.582}\left(\dfrac{(27.54-3\times1.581)^3}{96}-\dfrac{(3\times1.581-4.23)^3}{1.582\times0.582}\right)=50.97\cdots<51$, as desired.

In case 3), for $\gamma=3.51$ we show
$\dfrac{(3\lambda-\gamma)^3}{\lambda(\lambda-1)(\lambda-3)}-\dfrac{(27-2\gamma-3\lambda)^3}{48(\lambda-3)}<51$
for $\frac \gamma3<\lambda\le9-\frac23\gamma$ and
$\dfrac{(3\lambda-\gamma)^3}{\lambda(\lambda-1)(\lambda-3)}<51$ for
$\lambda\ge9-\frac23\gamma$.
The latter is easy.
In the former case we apply (3.11) for $\alpha=1$, $\beta=3$, $\gamma=3.51$.
We have
$\varphi(\lambda)=9\lambda^4-83.88\lambda^3+63.6804\lambda^2+112.8816\lambda-98.5608$ and
$\varphi(1.253)>0>\varphi(1.254)$, so
$f(\lambda)<\dfrac1{3-1.254}\left(\dfrac{(19.98-3\times1.253)^3}{48}-\dfrac{(3\times1.253-3.51)^3}{1.254\times0.254}\right)=50.899\cdots<51$.

(3.13) Now we divide the cases according to the nature of the place $P$.\nl
1) A representative pair of $P$ is obtained by a process as (3.3;1) and
$P=\SP(E_{(i)},Y_{(i)})$.\nl
2) The representing process goes via (3.3;2), so the locus of $P$ on
$E_{(i.j)}$ is contained in $Y_{(i.j)}$.

In the former case 1), we have $c_{(i)}\le1$ and
$m_{(1)}+\cdots+m_{(i)}=\mu_{(i)}\ge r_{(i)}+1=i+1$ (cf. (3.6)).
By (3.10), we can get rid of such a situation when $i\le5$, similarly as in
(3.1) and (3.2).
When $i>5$, we have $\mu_{(5)}=m_{(1)}+\cdots+m_{(5)}\ge\frac5i\mu_{(i)}>5$
since $m_{(1)}\ge\cdots\ge m_{(i)}$, hence (3.10;5) takes care of this case.

Thus, from now on, we consider the above case 2).

(3.14) Suppose that $i>5$.
By (3.7) we have
$1\le\Delta_{(i.1)}=\delta_{(i.1)}+m_{(i-1)}=\mu_{(i.1)}-r_{(i.1)}+m_{(i-1)}=m_{(1)}+\cdots+m_{(i-1)}+m_{(i.1)}-i+m_{(i-1)}$.
Hence
$\mu_{(5)}=m_{(1)}+\cdots+m_{(5)}\ge\frac5{i+1}(m_{(1)}+\cdots+m_{(i-1)}+m_{(i.1)}+m_{(i-1)})\ge5$, so this case is treated by (3.10;5).
Note that, if we choose a Hironaka model for which $mu_{(5)}<5$, then the above
situation cannot occur at {\it every} place whose locus is contained in
$Y_{(5)}$.

{}From now on, we assume $i\le5$.

(3.15) There are now the following three cases:\nl
1) $P=\SP(Y_{(i.j)})$.\nl
2) The locus of $P$ is the point $Y_{(i.j)}\cap Y_{(i.j-1)}$.
The process continues via a step (3.3;3).\nl
3) The locus of $P$ is a point of $Y_{(i.j)}$ off the othe $(-1)$-curves.
The process continues via a step as (3.3;5).

We have $\delta_{(i.j)}\ge1$ in case 1), $\delta_{(i.j)}+m_{(i.j-1)}\ge1$ in
case 2) by (3.7), and $\delta_{(i.j)}+m_{(i.j)}\ge1$ in case 3).
In any case $\delta_{(i.j)}+m_{(i.j-1)}\ge1$.

(3.16) First suppose that $i=5$.
We have
$1\le\delta_{(5.j)}+m_{(5.j-1)}=j\delta_{(4)}+m_{(5.1)}+\cdots\almb+m_{(5.j)}-j+m_{(5.j-1)}$.
On the other hand
$\mu_{(4)}=m_{(1)}+\cdots+m_{(4)}\ge4m_{(4)}\ge4(m_{(5.1)}+\cdots+m_{(5.j)})\ge4(j-1)m_{(5.j-1)}$, so
$1\le j(\mu_{(4)}-4)+\frac14\mu_{(4)}-j+\frac1{4(j-1)}\mu_{(4)}$, hence
$\mu_{(4)}\ge(5j+1)/(j+\frac14+\frac1{4(j-1)})=4(5j+1)(j-1)/j(4j-3)$.
This is an increasing function in $j$, so $\mu_{(4)}\ge\frac{22}5=4.4$ for
$j\ge2$.
By (3.10;4) we are done in this case.

Next suppose that $i=4$.
As above, we have
$1\le\delta_{(4.j)}+m_{(4.j-1)}=j\delta_{(3)}+m_{(4.1)}+\cdots\almb+m_{(4.j)}-j+m_{(4.j-1)}$ and
$\mu_{(3)}\ge3m_{(3)}\ge3(m_{(4.1)}+\cdots+m_{(4.j)})\ge3(j-1)m_{(4.j-1)}$.
Hence
$\mu_{(3)}\ge(4j+1)/(j+\frac13+\frac1{3(j-1)})=3(4j+1)(j-1)/j(3j-2)\ge\frac{26}7=3.7\cdots$ for $j\ge3$, and we are done in this case by (3.10;3).
But the case $j=2$ still survives.

Suppose now that $i=3$.
In this case we have
$\mu_{(2)}\ge(3j+1)/(j+\frac12+\frac1{2(j-1)})=2(3j+1)(j-1)/j(2j-1)$ as above.
Hence $\mu_{(2)}\ge\frac{128}{45}=2.84\cdots$ when $j\ge5$, and we can use
(3.10;2).
But the cases $j\le4$ are left to be studied.

When $i=2$, the same argument works if $j>108$.
Thus, we may assume $i\le4$, and further $j=2$ (resp. $j\le4$, $j\le108$) when
$i=4$ (resp. $i=3$, $i=2$).

(3.17) Suppose that $i=4$ and $j=2$.
We divide this case as in (3.15).

In case 1), we have $1\le\delta_{(4.2)}=2\delta_{(3)}+m_{(4.1)}+m_{(4.2)}-2$,
so $\mu_{(3)}\ge9/(2+\frac13)=3.85\cdots$ and we are done by (3.10;3).

In case 3), we have $1\le\delta_{(4.2)}+m_{(4.2)}$ by (3.7), and
$\mu_{(3)}\ge3m_{(3)}\ge3(m_{(4.1)}+m_{(4.2)})\ge6m_{(4.2)}$.
Hence $\mu_{(3)}\ge9/(2+\frac13+\frac16)=3.6$ and we are done.

In case 2), we have $1\le\Delta_{(4.2.k)}=\delta_{(4.2.k)}+m_{(4.2.k-1)}$.
By the method (3.6), we get
$\delta_{(4.2.k)}=\delta_{(4.2.1)}+(k-1)\delta_{(4.1)}+m_{(4.2.2)}+\cdots+m_{(4.2.k)}-(k-1)
=\delta_{(3)}+k\delta_{(4.1)}+m_{(4.2.1)}+\cdots+m_{(4.2.k)}-k
\le\delta_{(3)}+k\delta_{(4.1)}+m_{(4.1)}-k
=(k+1)\mu_{(4.1)}-5k-3$ and
$\mu_{(4.1)}\ge4m_{(4.1)}\ge4(m_{(4.2.1)}+\cdots+m_{(4.2.k)})\ge4(k-1)m_{(4.2.k-1)}$.
Therefore
$\mu_{(4.1)}\ge(5k+4)/(k+1+\frac1{4(k-1)})=4(5k+4)(k-1)/(4k^2-3)\ge\frac{56}{13}=4.3\cdots$ since $k\ge2$.
Hence we can use (3.10;4).

Thus we are done in this case $(i,j)=(4,2)$.

(3.18) It turns out that (3.10) is not enough to study the remaining cases with
$i\le3$.
We need a few results on $E_{(i.j)}$ of the type (3.10).
The counterpart of (3.9) in this situation is the following

\Lem {\sl Let $E_{(i.j)}$ be the surface obtained by the process (3.3;2), $H$
be the pull-back of $\SO_E(1)$, $Y_n$ and $Y_{i.n}$ be the total transforms of
the $(-1)$-curves $Y_{(n)}$ and $Y_{(i.n)}$ respectively, while their proper
transformations are denoted by $\tilde Y_n$ and $\tilde Y_{i.n}$.
Put
$L(u,w_1,\cdots,w_{i-1},w_{i.1},\cdots,w_{i.j})\almb=uH-w_1Y_1-\cdots-w_{i-1}Y_{i-1}-w_{i.1}Y_{i.1}-\cdots-w_{i.j}Y_{i.j}\in\Pic(E_{(i.j)})$,
$h^p(u,w_1,\cdots)=\dm H^p(E_{(i.j)},\almb L(u,w_1,\cdots))$,
$w=j(w_1+\cdots+w_{i-1})+w_{i.1}+\cdots+w_{i.j}$ and $\beta=(i-1)j+1$.
Then\nl
1) $h^0(u,w_1,\cdots)=0$ if $w>\beta u$.\nl
2) $h^1(u,w_1,\cdots)=0$ and
$h^0(u,w_1,\cdots)=\frac12(u+1)(u+2)-\frac12\sum_*w_*(w_*+1)$ if $u\ge
w_1+\cdots\almb+w_{i-1}+w_{i.1}$, $w_1\ge\cdots\ge w_{i-1}\ge
w_{1.1}+\cdots+w_{i.j}$ and $w_{i.1}\ge\cdots\ge w_{i.j}$.\nl
3) $h^0(u,w_1,\cdots)\le\frac12(1+u)(2+u)-\frac1{2j\beta
}(w+\frac{ij}2)^2+\frac18(i+j-1)$ if $w\le jw$.\nl
4) $h^0(u,w_1,\cdots)\le q(\beta u-w)$ if $w\le \beta u$, where\nl
$q(\theta)=_{\text{def}}\dfrac1{2\beta (\beta -j)}\theta^2+\dfrac{2\beta
-j+1}{2\beta (\beta -j)}\theta+1+\dfrac{(i-1)(j-1)^2}{8\beta }$.}

\Prf We argue as in (3.9) by induction on $u$.
The assertions are obvious when $u=0$, so we may suppose $u>0$.
Let $\tilde Z\in\vert H-Y_1-\cdots-Y_a\vert$ be the proper transform of a line
as in (3.9), where $2\le a\le i$.
The problem is reduced to the case $w_1\ge\cdots\ge w_{i-1}$ and
$w_{i.1}\ge\cdots\ge w_{i.j}$.
Moreover we may assume $w_{i-1}\ge w_{i.1}+\cdots+w_{i.j}$, since otherwise we
have $L(u, w_1,\cdots)\tilde Y_{i-1}<0$ and $h^0(u,w_1,\cdots)\le
h^0(u,w_1,\cdots,w_{i-1}+1,w_{i.1}+1,\cdots,w_{i.j}-1)$.

As for 1), we have $L\tilde Z=u-w_1-\cdots-w_a\le u-\frac w{ij}a\le\frac
u{ij}(ij-a\beta )<0$, so we are done as in (3.9).

2) is valid under the weaker assumption $u\ge w_1+\cdots+w_a$ instead of $u\ge
w_1+\cdots+w_{i-1}+w_{i.1}$.
We can prove this stronger assertion by induction on $u$.
But the cases $a=i$ and $a<i$ need slightly different arguments, so suppose
first $a<i$.
Any way we have $h^1(\tilde Z,L(u,\cdots))=0$, hence $h^1(u,\cdots)\le
h^1(L(u,\cdots)-\tilde Z)=h^1(u-1,w_1-1,\cdots,w_a-1,w_{a+1},\cdots)$.
This last term vanishes if $w_a>w_{a+1}$ by the induction hypothesis.
If not, we have $w_a=w_{a+1}$ and $h^1(\tilde Y_a,L(u,\cdots)-\tilde Z)=0$, so
the above term $\le h^1(u-1,\cdots,w_a,w_{a+1}-1,\cdots)$ as in (3.9).
Proceeding similarly we further get $\le h^1(u-1,w'_1,\cdots,w'_{a+1},\cdots)$
such that $w'_1\ge\cdots\ge w'_{a+1}=w_{a+1}-1$.
If $w_{a+1}>w_{a+2}$ we are done.
Otherwise by similar process we get $\le h^1(u-1,\almb
w''_1,\cdots,w''_{a+2},\cdots)$ such that $w''_1\ge\cdots\ge
w''_{a+2}=w_{a+2}-1$.
After several steps we eventually get $\le
h^1(u-1,\cdots,w_{i-1}-1,w_{i.1},\cdots)$.
If $w_{i-1}>w_{i.1}+\cdots+w_{i.j}$ then we are done.
If not, we have $w_{i-1}=w_{i.1}+\cdots+w_{i.j}$ and $h^1(\tilde
Y_{i-1},L(u-1,\cdots,w_{i-1}-1,w_{i.1},\cdots))=0$, so the above term $\le
h^1(L(\cdots,w_{i-1}-1,w_{i.1},\cdots)-\tilde
Y_{i-1})=h^1(u-1,\cdots,w_{i-1},w_{i.1}-1,\cdots,w_{i.j}-1)$.
This vanishes by the induction hypothesis.

Now we consider the case $a=i$.
In this case we have $h^1(u,\cdots)\le
h^1(u-1,w_1-1,\almb\cdots,w_{i-1}-1,w_{i.1}-1,w_{i.2},\cdots)$.
This vanishes by the induction hypothesis if $w_{i.1}>w_{i.2}$.
If $w_{i.1}=w_{1.2}$, we have $h^1(\tilde
Y_{i.1},L(u-1,\cdots,w_{i.1}-1,w_{i.2},\cdots))=0$, so the above $h^1\le
h^1(L(\cdots,w_{i.1}-1,w_{i.2},\cdots)-\tilde
Y_{i.1})=h^1(\cdots,w_{i.1},w_{1.2}-1,\cdots)$.
This vanishes unless $w_{i.2}=w_{i.3}$.
Continueing similarly, we eventually get $\le h^1(\cdots,w_{i.j-1},w_{i.j}-1)$,
which vanishes by the induction hypothesis.

Thus $h^1=0$ is proved in either case.
The assertion for $h^0$ follows from the Riemann-Roch theorem.

3) follows from 2) by the Schwarz inequality
$(w+\frac{ij}2)^2\le((w_1+\frac12)^2+\cdots+(w_{i-1}+\frac12)^2\almb+(w_{i.1}+\frac12)^2+\cdots+(w_{i.j}+\frac12)^2)(j^2+\cdots+j^2+1^2+\cdots+1^2)$.

4) is proved as in (3.9).
First consider the case $i=2$.
If $w\le ju$ 3) implies 4) since $q(\beta u-w)-(\text{the right side of
3})=\frac1{2j(\beta -j)}(w-ju)(w-ju+(i-3)j)\ge0$.
If $w>ju$ we have $u<w_1+w_{2.1}$, so $h^0(\tilde Z,L(u,\cdots))=0$ and
$h^0(u,\cdots)=h^0(u-1,w_1-1,w_{2.1}-1\cdots)\le q(\beta
(u-1)-j(w_1-1)-(w_{2.1}-1)-\cdots)=q(\beta u-w)$ by the induction hypothesis.
Thus we are done in either case.

Second consider the case $i=3$.
We have $q(\beta u-w)\ge(\text{the right side of 3})$ again, so the assertion
is true if $u\ge w_1+\cdots+w_a$.
If $u<w_1+\cdots+w_a$, we have
$h^0(u,\cdots)=h^0(u-1,w_1-1,\cdots,w_a-1,w_{a+1},\cdots)\le q(\beta
(u-1)-j(w_1-1)-\cdots)$ by the induction hypothesis.
When $a<i=3$, we have $a=2$ and $\beta (u-1)-j(w_1-1)-\cdots=\beta u-w-\beta
+aj=\beta u-w-1$, so we are done since $q(\beta u-w-1)-q(\beta u-w)=(-2(\beta
u-w)-2\beta +j)/2\beta (\beta -j)<0$.
When $a=i=3$, $\beta (u-1)-j(w_1-1)-j(w_2-1)-(w_{3.1}-1)-w_{3.2}-\cdots=\beta
u-w$ and we are done.

\Rmk 1),2) and 3) are true even if $i>3$, but 4) is uncertain.

(3.19) \Lem {\sl $\pi_1^*B-\frac92E$ is big in any of the following cases:\nl
1) $\Xi_{(3.2)}\ge6.31$.\nl
2) $\Xi_{(3.3)}\ge9.13$.\nl
3) $\Xi_{(3.4)}\ge12$.}

\Prf Similar as (3.10).
For $\gamma=6.31$ (resp. $9.13$, $12$) in case $j=2$ (resp. $j=3$, $j=4$), we
show
$$\int_{\eta_\lambda(x)\le
jx}(\frac12x^2-\frac{\eta_\lambda(x)^2}{2j\beta})dx+\int_{jx\le\eta_\lambda(x)\le \beta x}\frac{\beta x-\eta_\lambda(x))^2}{2\beta (\beta -j)}dx<\frac{51}6 ,$$
where $\beta =2j+1$ is as in (3.18) and
$\eta_\lambda(x)=\Max(0,\lambda(x-3)+\gamma)$.
The problem is reduced to $f_{j,\beta ,\gamma}(\lambda)<51$ if $\frac
\gamma3<\lambda\le3\beta -\frac23\gamma$, and
$(3\lambda-\gamma)/\lambda(\lambda-j)(\lambda-\beta )<51$ if $\lambda\ge3\beta
-\frac23\gamma$.
The latter inequality is easy to prove.
To show the former, we argue as in (3.12) using (3.11).
$\varphi(t)$ changes its sign from $+$ to $-$ once at some $\lambda_1$, and
$f_{j,\beta ,\gamma}(t)$ attains the maximum at $\lambda_1$.
By computation we see $2.161<\lambda_1<2.162$ (resp. $3.069<\lambda_1<3.07$) in
case $j=2$ (resp. $j=3$), and get the desired inequality as before.
The case $j=4$, $\gamma=12$ is even easier.
Indeed $\varphi<0$ and the maximum is attained at $\lambda_1=4$.

(3.20) \Lem {\sl Suppose that $\delta_0<1$ and $m_0<1$ in the situation
(3.3;5). Then\nl
1) $\delta_0+m_0\ge\frac32$.\nl
2) $\delta_0+\dfrac{q^2+q+1}{(q+1)^2}m_0\ge2-\dfrac1{q+1}$ for some $q\ge2$.\nl
3) $m_0\ge\frac67$.}

\Prf The blow-up process continues via a step
$$\vdash-[2]-\cdots-[i'-1]-[i'.p]-\cdots-[i'.1]$$
unless $P=\SP(Y_{[i']})$, where we set $p=j'$.
Since $1>m_0\ge m_{[2]}\ge\cdots\ge m_{[i']}$ (cf. (3.5), (3.6)), we have
$\delta_0>\delta_0+m_{[2]}-1=\delta_{[2]}>\cdots>\delta_{[i']}$, so the latter
exceptional case is ruled out.
Moreover, replacing $\delta_0$ and $m_0$ by $\delta_{[i'-1]}$ and $m_{[i'-1]}$
if necessary, we can suppose $i'=2$.
Then we have
$1\le\delta_{[2.p]}+m_{[2.p-1]}=p(\delta_0-1)+m_{[2.1]}+\cdots+m_{[2.p]}+m_{[2.p-1]}$ by (3.7).
Since $1>m_0\ge m_{[2.1]}+\cdots+m_{[2.p]}\ge (p-1)m_{[2.p-1]}$,
we infer $1\le p(\delta_0-1)+(1+\frac1{p-1})m_0\le 2(\delta_0-1)+2m_0$,
which yields 1).

If $\delta_{[2.p]}\ge1$, we have $1\le
p(\delta_0-1)+m_{[2.1]}+\cdots+m_{[2.p]}\le p(\delta_0-1)+m_0$, contradicting
the assumption.

If $\delta_{[2.p]}+m_{[2.p]}\ge\frac32$ (cf. 1), then $\frac32\le
p(\delta_0-1)+m_0+\frac1{2p}m_0<1+\frac1{2p}$, absurd.

Thus the process continues via the step (3.3;6), where we set $k'=q$. Then
$1\le\delta_{[2.p.q]}+m_{[2.p.q-1]}\le\delta_0+q(\delta_{[2.p-1]}-1)+m_{[2.p.1]}+\cdots+m_{[2.p.q]}+m_{[2.p.q-1]}$ and\nl
$\delta_{[2.p-1]}\le(p-1)(\delta_0-1)+m_{[2.1]}+\cdots+m_{[2.p-1]}$.
Hence\nl
$1\le((p-1)q+1)\delta_0-pq+q(m_{[2.1]}+\cdots+m_{[2.p-1]})+m_{[2.p.1]}+\cdots+m_{[2.p.q]}+m_{[2.p.q-1]}$.
On the other hand $m_0\ge
m_{[2.1]}+\cdots+m_{[2.p-1]}+m_{[2.p.1]}\ge(p-1)m_{[2.p-1]}+m_{[2.p.1]}\ge((p-1)q+1)m_{[2.p.q]}$ and
$qm_{[2.p.1]}\ge m_{[2.p.1]}+\cdots+m_{[2.p.q-2]}+2m_{[2.p.q-1]}$, so we get\nl
$1+pq\le((p-1)q+1)\delta_0+(q+\frac1{(p-1)q+1})m_0\le(q+1)\delta_0+(p-2)q+(q+\frac1{q+1})m_0$.
Hence $1+2q\le(q+1)\delta_0+(q+\frac1{q+1})m_0$, which yields 2).

2) yields $dfrac{(q^2+q+1)}{(q+1)^2}m_0\ge1-\dfrac1{q+1}=dfrac q{q+1}$, so
$m\ge\dfrac{q(q+1)}{q^2+q+1}=1-\dfrac1{q^2+q+1}\ge\dfrac67$, proving 3).

(3.21) Let us consider the case $i=3$, $j\le 4$.

If $P=\SP(Y_{(3.j)})$, we have $1\le\delta_{(3.j)}$ and $\mu_{(3.j)}\ge
r_{(3.j)}+1=3j+1$, so (3.19) applies.

In case (3.3;5), we have $\delta_{(3.j)}+m_{(3.j)}\ge\frac32$ by (3.20;1).
Since $\delta_{(3.j)}=\mu_{(3.j)}-3j$ and
$\mu_{(3.j)}=j(m_{(1)}+m_{(2)})+m_{(3.1)}+\cdots+m_{(3.j)}\ge 2j
m_{(2)}+m_{(3.1)}+\cdots+m_{(3.j)}\ge(2j+1)(m_{(3.1)}+\cdots+m_{(3.j)})\ge(2j+1)jm_{(3.j)}$, we have $\mu_{(3.j)}\ge(3j+\frac32)/(1+\frac1{j(2j+1)})$.
Hence (3.19) applies in any case $j=2,3,4$.

In the remaining cases the process continues via (3.3;3).
Then we have
$1\le\delta_{(3.j.k)}+m_{(3.j.k-1)}=\delta_{(2)}+k(\delta_{(3.j-1)}-1)+m_{(3.j.1)}+\cdots+2m_{(3.j.k-1)}+m_{(3.j.k)}\le
\delta_{(2)}+k(\delta_{(3.j-1)}-1)\almb+km_{(3.j.1)}+m_{(3.j.k)}=
((j-1)k+1)(m_{(1)}+m_{(2)})+k(m_{(3.1)}+\cdots+m_{(3.j-1)}+m_{(3.j.1)})+m_{(3.j.k)}-2-3(j-1)k-k
=\frac{(j-1)k+1}j\mu_{(3.j)}+\frac{k-1}j(m_{(3.1)}+\cdots+m_{(3.j)})+m_{(3.j.k)}-(3j-2)k-2$ and
$\mu_{(3.j)}\ge(2j+1)(m_{(3.1)}+\cdots+m_{(3.j)})\ge(2j+1)((j-1)m_{(3.j-1)}+m_{(3.j.1)})\ge(2j+1)\almb((j-1)k+1)m_{(3.j.k)}$.
Hence
$\mu_{(3.j)}\ge((3j-2)k+3)/(\frac{(j-1)k+1}j+\frac{k-1}j\frac1{2j+1}+\frac1{(2j+1)((j-1)k+1)})$.
When $j=2$, this yields $\mu_{(3.j)}\ge5(4k+3)/(3k+2+\frac1{k+1})\ge 6.6$ for
every $k\ge2$.
When $j=3$ we get $\mu_{(3.j)}\ge7(7k+3)/(5k+2+\frac1{2k+1})>9.7$ for every
$k\ge2$.
When $j=4$ we get $\mu_{(3.j)}\ge9(10k+3)/(7k+2+\frac1{3k+1})>12$ for every
$k\ge2$.
Thus (3.19) applies in any case.

(3.22) Now it remains only the cases $i=2$, $j\le108$.
We will show the bigness of $\pi_1^*B-\frac92E$ assuming $\Xi_{(2.j)}\ge
\gamma_{(2.j)}$ for some number $\gamma=\gamma_{(2.j)}$ slightly smaller than
$r_{(2.j)}+1=2j+1$.

As before,
$\Xi_{(2.j)}(\pi_1^*B-xE)\ge\eta_\lambda(x)=_{\text{def}}\Max(0,\lambda(x-3)+\gamma)$ for some $\lambda\ge\frac \gamma3$.
By (3.18), we should estimate
$$\psi(\lambda)=\int_{\eta_\lambda(x)\le
jx}\left(\frac12x^2-\frac{\eta_\lambda(x)^2}{2j\beta}\right)dx+\int_{jx\le\eta_\lambda(x)\le \beta x}\frac{\beta x-\eta_\lambda(x))^2}{2\beta (\beta -j)}dx\ .$$
Unlike the case $i=3$, we have $\frac \gamma3<3j-\frac23\gamma$, so there are
the three cases:\nl
1) $\frac \gamma3<\lambda\le3j-\frac23\gamma$.
$6\psi(\lambda)={\dsize\int}_0^{3-\frac
\gamma\lambda}3x^2+{\dsize\int}_{3-\frac
\gamma\lambda}^\frac92\left(3x^2-\dfrac{3(\lambda(x-3)+\gamma)^2}{j\beta}\right)dx$\nl
$=\dfrac{729}8-\dfrac1{j\beta
\lambda}\left(\dfrac32\lambda+\gamma\right)^3$.\nl
2) $3j-\frac23\gamma\le\lambda\le3\beta -\frac23\gamma$.
$6\psi(\lambda)={\dsize\int}_0^{3-\frac
\gamma\lambda}3x^2dx+{\dsize\int}_{3-\frac
\gamma\lambda}^{\frac{3\lambda-\gamma}{\lambda-j}}\left(3x^2-\dfrac{3\eta_\lambda(x)^2}{j\beta }\right)dx+$\nl
${\dsize\int}_{\frac{3\lambda-\gamma}{\lambda-j}}^\frac92\dfrac{3(\beta
x-\eta_\lambda(x))^2}{(\beta -j)\beta }dx
=\dfrac{(3\lambda-\gamma)^3}{\lambda(\lambda-j)(\lambda-\beta )}-\dfrac{(9\beta
-2\gamma-3\lambda)^3}{8\beta (\beta -j)(\lambda-\beta )}$.\nl
3) $3\beta -\frac23\gamma\le\lambda$.
$6\psi(\lambda)=\dfrac{(3\lambda-\gamma)^3}{\lambda(\lambda-j)(\lambda-\beta
)}$.

In the range 3), we have
$\frac
d{d\lambda}\log\psi=\frac9{3\lambda-\gamma}-\frac1\lambda-\frac1{\lambda-j}-\frac1{\lambda-\beta }<0$,
since $\frac1\lambda+\frac1{\lambda-j}+\frac1{\lambda-\beta
}>3\root3\of{\frac1\lambda\cdot\frac1{\lambda-j}\cdot\frac1{\lambda-\beta }}$
and $\root3\of{\lambda(\lambda-j)(\lambda-\beta
)}\le\frac13(\lambda+(\lambda-j)+(\lambda-\beta
))=\frac13(3\lambda-2j-1)\le\frac13(3\lambda-\gamma)$.

In the range 1), $\frac
d{d\lambda}\log(\frac1\lambda(\frac32\lambda+\gamma)^3)=\frac9{(3\lambda+2\gamma)}-\frac1\lambda=2(3\lambda-\gamma)/(3\lambda+2\gamma)\lambda\ge0$, so $\psi(\lambda)$ decreases in this range too.

In the range 2), we use (3.11) to analyze $6\psi=f_{j,\beta ,\gamma}$.
Note that $\beta =j+1$ and $\epsilon=2j+1-\gamma\ge0$.
Set $c=3j-\frac23\gamma$ and $d=3\beta -\frac23\gamma$.
We have $\varphi^{'''}\ge\varphi^{'''}(c)=36(5j-8+2\epsilon)>0$, so $\varphi''$
increases.
Since $\varphi''(d)=50j^2+236j+236+(40j+80)\epsilon+8\epsilon^2>0$, $\varphi''$
is either everywhere positive in this range 2) or varies from $-$ to $+$.
In either case, since $\varphi'(c)=-12(11j-2+2\epsilon)(j-1+\epsilon)<0$ and
$\varphi'(d)=18(j+1)(j+2-2\epsilon)>0$, $\varphi'$ varies from $-$ to $+$.
We have $\varphi(c)=-36j(j-1+\epsilon)^2<0$ and
$\varphi(d)=-9(j+1)(4j^2+7j+7+8(j+2)\epsilon+4\epsilon^2)<0$, so $\varphi<0$ in
this range and $\psi$ decreases.

Thus we conclude $\psi(\lambda)\le\psi(\frac
\gamma3)=\frac{27}{16}(9-\frac{\gamma^2}{j\beta })$, hence

(3.23) {\sl $\pi_1^*B-\frac92E$ is big if $\Xi_{(2.j)}>\sqrt{321j(j+1)}/9$.}

Indeed, $\frac{27}{16}(9-\frac{\gamma^2}{j\beta })<\frac{51}6$ iff
$\gamma>\sqrt{321j\beta }/9$.

(3.24) By (3.7) we have
$1\le\delta_{(2.j)}+m_{(2.j-1)}=\mu_{(2.j)}-2j+m_{(2.j-1)}$, while
$\mu_{(2.j)}=jm_{(1)}+m_{(2.1)}+\cdots+m_{(2.j)}\ge(j+1)(j-1)m_{(2.j-1)}$.
Hence $\mu_{(2.j)}\ge(2j+1)/(1+\frac1{(j+1)(j-1)})$.
By computation we see $(\text{the right side})\ge\sqrt{321j(j+1)}/9$ for
$j\ge14$.
Now the cases $j\le13$ are left.

(3.25) We assume $j\le13$ from now on.
If $\delta_{(2.j)}\ge1$, then $\mu_{(2.j)}\ge2j+1$ and we are done by (3.23)
for any $j\ge2$.
If the process continues as (3.3;5), we have
$\delta_{(2.j)}+(1-\frac q{(q+1)^2})m_{(2.j)}\ge2-\frac1{q+1}$ for some $q\ge2$
by (3.20;2).
Moreover $\mu_{(2.j)}\ge(j+1)(m_{(2.1)}+\cdots+m_{(2.j)})\ge j(j+1)m_{(2.j)}$,
so
$\left(1+\dfrac{q^2+q+1}{(q+1)^2j(j+1)}\right)\mu_{(2.j)}\ge2j+2-\frac1{q+1}$.
We may assume $\mu_{(2.j)}<2j+1$, hence
$\dfrac{q^2+q+1}{(q+1)^2}\cdot\dfrac{2j+1}{j(j+1)}\ge 1-\frac1{q+1}=\frac
q{q+1}$, so $\dfrac{2j+1}{j(j+1)}\ge\dfrac{q(q+1)}{q^2+q+1}$.
But the left side $\le\frac56$ for any $j\ge2$, while the right side
$=1-\frac1{q^2+q+1}\ge\frac67$ for any $q\ge2$.
Thus we get a contradiction.

Therefore, the blow-up process must continue as (3.3;3).

(3.26) Now we have $1\le\delta_{(2.j.k)}+m_{(2.j.k-1)}$.
Since
$\delta_{(2.j.k)}=\delta_{(1)}+k(\delta_{(2.j-1)}-1)+m_{(2.j.1)}+\cdots+m_{(2.j.k)}$,
$\delta_{(2.j-1)}=(j-1)(\delta_{(1)}-1)+m_{(2.1)}+\cdots+m_{(2.j-1)}$,
$\delta_{(1)}=m_{(1)}-1$ and
$\mu_{(2.j)}=jm_{(1)}+m_{(2.1)}+\cdots+m_{(2.j-1)}+m_{(2.j.1)}\ge(j+1)(m_{(2.1)}+\cdots+m_{(2.j-1)}+m_{(2.j.1)})\ge(j+1)((j-1)k+1)m_{(2.j.k)}$,
we infer
$(2j-1)k+2\le
((j-1)k+1)m_{(1)}+k(m_{(2.1)}+\cdots\almb+m_{(2.j-1)})+m_{(2.j.1)}+\cdots+2m_{(2.j.k-1)}+m_{(2.j.k)}\le
\frac{(j-1)k+1}j\mu_{(2.j)}+\frac{k-1}j(m_{(2.1)}+\cdots+m_{(2.j-1)})-\frac{(j-1)(k-1)}j m_{(2.j.1)}+m_{(2.j.2)}+\cdots\le
(\frac{(j-1)k+1}j+\frac{k-1}{j(j+1)})\mu_{(2.j)}-(k-1)m_{(2.j.1)}+m_{(2.j.2)}+\cdots\almb+2m_{(2.j.k-1)}+m_{(2.j.k)}\le
\frac{jk+1}{j+1}\mu_{(2.j)}+m_{(2.j.k)}\le
(\frac{jk+1}{j+1}+\frac1{(j+1)((j-1)k+1)})\mu_{(2.j)}$, therefore
$\mu_{(2.j)}\ge((2j-1)k+2)/(\frac{jk+1}{j+1}+\frac1{(j+1)((j-1)k+1)})=
(j+1)((2j-1)k+2)((j-1)k+1)/((jk+1)((j-1)k+1)+1)$.
Hence $\mu_{(2.j)}\ge(j+1)(2j-1)/j$, since
$j((2j-1)k+2)((j-1)k+1)-(2j-1)((jk+1)((j-1)k+1)+1)=(j-1)(k-2)\ge0$.

(3.27) Now we divide the cases by $j$.

For $j=13$, we have
$\mu_{(2.j)}\ge14(25k+2)(12k+1)/(156k^2+25k+2)\ge\frac{350}{13}=26.9\cdots$ for
any $k\ge2$, so we are done by (3.23) since $\sqrt{321j(j+1)}/9=26.85\cdots$.

For $j=12$, we have
$\mu_{(2.j)}\ge13(23k+2)(11k+1)/(132k^2+23k+2)\ge\frac{299}{12}=24.9\cdots$ for
any $k\ge2$.
For $j=11$,
$\mu_{(2.j)}\ge12(21k+2)(10k+1)/(110k^2+21k+2)\ge\frac{252}{11}=22.9\cdots$.
For $j=10$, $\mu_{(2.j)}\ge20.9$ and for $j=9$
$\mu_{(2.j)}\ge\frac{170}9=18.888\cdots$.
In any of these cases we can apply (3.23).

However, when $j\le8$, we must still work harder.
We need a result of the type (3.23) in terms of $\Xi_{(2.j.k)}$.

(3.28) \Lem {\sl Let
$L(u,w_1,w_{2.1},\cdots,w_{2.j-1},w_{2.j.1},\cdots,w_{2.j.k})=uH-\sum_*w_*Y_*\in\Pic(E_{(2.j.k)})$ as in (3.9) and (3.18).
Put $h^p(u,w_1,\cdots)=\dm H^p(E_{(2.j.k)},L(u,w_1,\cdots))$,
$w=((j-1)k+1)w_1+k(w_{2.1}+\cdots+w_{2.j-1})+w_{2.j.1}+\cdots+w_{2.j.k}$,
$\beta =jk+1$ and $\alpha=\beta-k$.
Then\nl
1) $h^0(u,\cdots)=0$ if $\beta u<w$.\nl
2) $h^1(u,\cdots)=0$ and
$h^0(u,\cdots)=\frac12(u+1)(u+2)-\sum_*\frac12w_*(w_*+1)$ if $u\ge
w_1+w_{2.1}$, $w_1\ge w_{2.1}+\cdots+w_{2.j-1}+w_{2.j.1}$, $w_{2.1}\ge\cdots\ge
w_{2.j-1}\ge w_{2.j.1}+\cdots+w_{2.j.k}$ and if $w_{2.j.1}\ge\cdots\ge
w_{2.j.k}$.\nl
3)
$h^0(u,\cdots)\le\dfrac12(u+1)(u+2)-\dfrac1{2\alpha\beta}(w+\dfrac{(2j-1)k+1}2)^2+\dfrac{j+k}8$ if $w\le\alpha u$.\nl
4) $h^0(u,\cdots)\le q(\beta u-w)$ if $w\le \beta u$, where
$q(\theta)=_{\text{def}}\dfrac1{2\beta k}\theta^2+\dfrac{\beta+k+1}{2\beta
k}\theta+1+\dfrac{j+k}8$.}

\Prf The same method works as in (3.18), since we may assume $w_1\ge
w_{2.1}+\cdots+w_{2.j-1}+w_{2.j.1}$, $w_{2.1}\ge\cdots\ge w_{2.j-1}\ge
w_{2.j.1}+\cdots+w_{2.j.k}$ and $w_{2.j.1}\ge\cdots\ge w_{2.j.k}$.
If $u\ge w_1+w_{2.1}$ further, then $\alpha u-w\ge\alpha
w_{2.1}-k(w_{2.1}+\cdots+w_{2.j-1})-w_{2.j.1}-\cdots-w_{2.j.k}\ge
w_{2.1}-w_{2.j.1}-\cdots-w_{2.j.k}\ge0$.
Details are left to the reader.

(3.29) Suppose that $\Xi_{(2.j.k)}\ge \gamma $ for some $\gamma =\gamma
_{(2.j.k)}\le2\beta -k$.
Then $\Xi_{(2.j.k)}(\pi_1^*B-xE)\ge\eta_\lambda(x)=\Max(0,\lambda(x-3)+\gamma
)$ for some $\lambda>\dfrac \gamma 3$.
We will estimate $\psi(\lambda)=\int_{\eta_\lambda(x)\le\alpha
x}(\frac12x^2-\frac{\eta_\lambda(x)^2}{2\beta \alpha })dx+\int_{\alpha
x\le\eta_\lambda(x)\le \beta x}\frac{(\beta x-\eta_\lambda(x))^2}{2\beta k}dx$.
Since $(3\alpha -\frac23\gamma )-\frac \gamma 3\ge \beta -2k=(j-2)k+1>0$, there
are the following three cases as in (3.22):\nl
1) $\frac \gamma 3<\lambda\le3\alpha -\frac23\gamma $.
$6\psi(\lambda)=\int_0^{3-\frac \gamma \lambda}3x^2dx+\int_{3-\frac \gamma
\lambda}^\frac92(3x^2-\frac{3\eta_\lambda(x)^2}{\alpha \beta })dx=
\frac{729}8-\frac1{\alpha \beta \lambda}(\frac32\lambda+\gamma )^3$.\nl
2) $3\alpha -\frac23\gamma\le\lambda\le 3\beta -\frac23\gamma $.
$6\psi(\lambda)=f_{\beta -k,\beta ,\gamma }(\lambda)$.\nl
3) $3\beta -\frac23\gamma \le\lambda$.
$6\psi(\lambda)=(3\lambda-\gamma )^3/\lambda(\lambda-\beta +k)(\lambda-\beta
)$.

As before, we easily checks that $\psi(\lambda)$ decreases in the range 1) and
3).
In the range 2), by (3.11) we get $\varphi^{'''}(c)=36(5\beta -13k+2\epsilon)$
and $\varphi^{'''}(d)=36(5\beta +5k+2\epsilon)>0$ for $c=3\alpha -\frac23\gamma
$, $d=3\beta -\frac23\gamma $ and $\epsilon=2\beta -k-\gamma \ge0$.
Hence $\varphi^{'''}$ is everywhere positive, or varies from $-$ to $+$.
We have $\varphi''(c)=50\beta ^2-404\beta k+482k^2+(40\beta
-176k)\epsilon+8\epsilon^2$ and
$\varphi''(d)=50\beta ^2136\beta k+50k^2+40(\beta +k)\epsilon+8\epsilon^2>0$.
Hence $\varphi''$ is everywhere positive, or varies from $-$ to $+$, or $+$
$\rightarrow$ $-$ $\rightarrow$ $+$.
This last case can occur only when $\varphi^{'''}(c)<0<\varphi''(c)$, so
$5\beta <13k$ and $j=2$, $\beta =2k+1$, but then
$\varphi''(c)=-126k^2-204k+50+(-96k+40)\epsilon+8\epsilon^2<0$ unless
$\epsilon$ is very large.
This case is thus ruled out if, e.g., $\epsilon<19$.
We have $\varphi'(c)=-12k(11\beta -13k+2\epsilon)(\beta -2k+\epsilon)<0$ and
$\varphi'(d)=18k\beta (\beta +k-2\epsilon)>0$ unless $2\epsilon>\beta
+k=(j+1)k+1$,
so usually $\varphi'$ varies from $-$ to $+$.
Since $\varphi(c)=-36k\alpha (\beta -2k+\epsilon)^2<0$ and
$\varphi(d)=-9\beta k((4\beta ^2-\beta k+4k^2)+8(\beta
+k)\epsilon+4\epsilon^2)<0$,
$\varphi<0$ in the range 2) if, e.g., $\epsilon\le3$.

Thus we usually obtain $\psi(\lambda)<\psi(\frac \gamma
3)=\frac{27}{16}(9-\frac{\gamma ^2}{\alpha \beta })$.
Hence

(3.30) {\sl $\pi_1^*B-\frac92E$ is big if $\Xi_{(2.j.n)}>\gamma_{(2.j.n)}$ and
if $(2j-1)n+2-\gamma_{(2.j.n)}\le3$, where
$\gamma_{(2.j.n)}=\sqrt{321((j-1)n+1)(jn+1)}/9$.}

(3.31) By (3.7) we have
$1\le\delta_{(2.j.k)}+m_{(2.j.k-1)}$, so
$(2j-1)k+2\le \mu_{(2.j.k)}+m_{(2.j.k-1)}$.
Since
$\mu_{(2.j.k)}=((j-1)k+1)m_{(1)}+k(m_{(2.1)}+\cdots+m_{(2.j-1)})+m_{(2.j.1)}+\cdots+m_{(2.j.k)}\ge
(jk+1)(m_{(2.1)}+\cdots+m_{(2.j-1)})+((j-1)k+2)m_{(2.j.1)}+m_{(2.j.2)}+\cdots+m_{(2.j.k)}\ge
(j-1)(jk+1)m_{(2.j-1)}+\cdots\ge
((j-1)(jk+1)(k-1)+(j-1)k+2+(k-2))m_{(2.j.k-1)}=
((j-1)jk^2-(j^2-3j+1)k+1-j)m_{(2.j.k-1)}$, we get
$$\mu_{(2.j.k)}\ge((2j-1)k+2)/\left(1+\frac1{(j-1)jk^2-(j^2-3j+1)k-j+1}\right).$$
As we will see later, this estimate is enough for (3.30) in most cases.
But here is a danger.

(3.32) Even if we show that the above lower bound $\ge\sqrt{321\alpha \beta
}/9$ for any $k$, this is not enough for our purpose !
Why so ?

We want to derive a contradiction assuming that $\pi_1^*B-\frac92E$ is not big.
Hence we may suppose $\Xi_{(2.j.k)}<\sqrt{321\alpha \beta }/9=\gamma$.
But this does not mean $\mu_{(2.j.k)}<\gamma$.
It means that if we replace the Hironaka model in (1.4) by another one of
possibly larger degree, then $\mu_{(2.j.k)}<\gamma$ for this model.
On this new model, $(2.j.k)$ cannot be an intermediate step of the blow-up
process representing a bad place as in (1.7).
However, there might exist a bad place for which $(2.j.k')$ is an intermediate
step for some $k'>k$.
This place can be killed by replacing the Hironaka model, but then there might
still other possibilities $(2.j.k'')$ with $k''>k'$ survive.
If such replacing processes repeat infinitely, we cannot complete the proof.

(3.33) In order to avoid the above danger, we will show that for each $j\le8$,
there is a number $n$ such that
$\delta_{(2.j.n)}\ge\sqrt{321((j-1)n+1)(jn+1)}/9$ if $k\ge n$.
Then, if we take a Hironaka model such that
$\delta_{(2.j.n)}<\gamma_{(2.j.n)}$, then the above possibilities in (3.32) are
killed for {\it all} $(2.j.k)$ with $k\ge n$ as (3.10;5) kills all the
possibility $(i)$ with $i\ge5$ (cf. (3.13)).

(3.34) As a matter of fact, $n=12$ has the above property.

To see this, we first note
$2\ge m_{(1)}\ge m_{(2.1)}+\cdots+m_{(2.j-1)}+m_{(2.j.1)}
\ge(j-1)m_{(2.j-1)}+m_{(2.j.1)}
\ge(j-1)(m_{(2.j.1)}+\cdots+m_{(2.j.k)})+m_{(2.j.1)}
\ge((j-1)(k-1)+1)m_{(2.j.k-1)}$, so
$1-\nomb\delta_{(2.j.k)}\le m_{(2.j.k-1)}\le 2/((j-1)(k-1)+1)\le\frac2j$
for any $k\ge2$.
Since $r_{(2.j.n)}=(2j-1)n+1$ for any $n$, we have
$\mu_{(2.j.k)}\ge(2j-1)k+2-\frac2j$.
We have
$\mu_{(2.j.n)}=((j-1)n+1)m_{(1)}+n(m_{(2.1)}+\cdots\almb+m_{(2.j-1)})+m_{(2.j.1)}+\cdots+m_{(2.j.n)}$ for any $n$, so
$\mu_{(2.j.n)}=\frac n k \mu_{(2.j.k)}+\frac{k- n}k m_{(1)}+\frac{k-
n}k(m_{(2.j.1)}+\cdots+m_{(2.j. n)}-\frac n k(m_{(2.j.
n+1)}+\cdots+m_{(2.j.k)})
\ge\frac n k\mu_{(2.j.k)}+\frac{k- n}k m_{(1)}$ if $ n\le k$, since
$m_{(2.j.1)}\ge\cdots\almb\ge m_{(2.j.k)}$ implies
$(m_{(2.j.1)}+\cdots+m_{(2.j. n)})/ n\ge(m_{(2.j. n+1)}+\cdots+m_{(2.j.k)})/(k-
n)$.
Therefore $\epsilon=_{\text{def}}1-\delta_{(2.j. n)}=(2j-1) n+2-\mu_{(2.j.
n)}\le
(2j-1) n+2-\frac n k((2j-1)k+2-\frac2j)-\frac{k- n}k m_{(1)}=\frac{k-
n}k(2-m_{(1)})+\frac{2 n}{jk}$.
In particular $\epsilon<2-m_{(1)}+\frac2j<3$, so it suffices to show
$(2j-1) n+2-\epsilon\ge\sqrt{321((j-1) n+1)(j n+1)}/9$ by (3.30).

{}From $1\le\delta_{(2.j.k)}+m_{(2.j.k-1)}$ we get
$(2j-1)k+2\le\mu_{(2.j.k)}+m_{(2.j.k-1)}
=((j-1)k+1)m_{(1)}+k(m_{(2.1)}+\cdots+m_{(2.j-1)})+m_{(2.j.1)}+\cdots+m_{(2.j.k)}+m_{(2.j.k-1)}
\le(jk+1)m_{(1)}-(k-1)m_{(2.j.1)}+m_{(2.j.2)}+\cdots+2m_{(2.j.k-1)}+m_{(2.j.k)}
\le(jk+1)m_{(1)}+m_{(2.j.k)}\le(jk+1+\frac1{(j-1)k+1})m_{(1)}$,
hence $m_{(1)}\ge((2j-1)k+2)/(jk+1+\frac1{(j-1)k+1})$ and
$2-m_{(1)}\le((j-1)k^2+k+2)/(jk+1)((j-1)k+1)\le\frac1j$ for any $k\ge3$.
Hence $\epsilon\le2-m_{(1)}+\frac2j\le\frac3j$.
For $ n=12$, this implies $2((2j-1) n+2)\epsilon<4j n\epsilon\le n^2$, which
yields
$(2j-1) n+2-\epsilon\ge2\sqrt{((j-1) n+1)(j n+1)}$.
This is more than enough.

(3.35) Now we consider the case $4\le k\le 12$.
By (3.31) we have $\epsilon=_{\text{def}}1-\delta_{(2.j.k)}$\nl
$=(2j-1)k+2-\mu_{(2.j.k)}\le((2j-1)k+2)/((j-1)jk^2-(j^2-3j+1)k-j+2)$\nl
$=((2j-1)k+2)/(jk+1)((j-1)(k-1)+1)$.

We claim $2((2j-1)k+2)^2\le k^2(jk+1)((j-1)(k-1)+1)$.
This yields $k^2-2((2j-\nomb 1)k+\nomb2)\epsilon\ge\nomb0$, hence
$\mu_{(2.j.k)}=(2j-1)k+2-\epsilon\ge2\sqrt{((j-1)k+1)(jk+1)}$,
which is more than enough.

To show the claim, we note $(j-1)(k-1)+1>\frac{j-1}j\cdot\frac{k-1}k(jk+1)$ and
$jk+1>\frac12((2j-1)k+2)$.
Hence it suffices to show $2\le k^2\frac14\frac{j-1}j\frac{k-1}k$, namely $8\le
k(k-1)\frac{j-1}j$.
Since $k\ge4$ and $j\ge2$, this is valid except $k=4$ and $j=2$.
But in this last case the claim is verified by direct computation.

Thus we are done in this case.

(3.36) Here let $k=3$.
By (3.31) we have $\mu_{(2.j.k)}\ge(6j-1)(6j^2-j-2)/(3j+1)(2j-1)$, and further
$>\sqrt{321(3j+1)(3j-2)}/9$ unless $j=2$.

Thus, in the remaining cases we have $k=2$ and $2\le j\le8$, or $(j,k)=(2,3)$.

(3.37) Suppose that $\delta_{(2.j.k)}\ge1$.
Then $\mu_{(2.j.k)}\ge(2j-1)k+2=((j-1)k+1)+(jk+1)\ge2\sqrt{((j-1)k+1)(jk+1)}$,
so (3.30) applies.
Hence the blow-up process does not stop.

If it continues as (3.3;5), then
$\mu_{(2.j.k)}\ge(j-1)(jk+1)m_{(2.j-1)}+((j-1)k+2)m_{(2.j.1)}+m_{(2.j.2)}+\cdots+m_{(2.j.k)}
\ge((j-1)(jk+1)k+((j-1)k+2)+k-1)m_{(2.j.k)}
=(jk+1)((j-1)k+\nomb1)m_{(2.j.k)}\ge\frac67(jk+1)((j-1)k+1)$,
which is more than enough to apply (3.30).

Thus the process must continue as (3.3;4).
In this case $1\le\delta_{(2.j.k.\ell)}+m_{(2.j.k.\ell-1)}$ by (3.7).

(3.38) Here we consider the case $j=2$, $k=3$.
Then $\mu_{(2.j.k)}=4m_{(1)}+3m_{(2.1)}+m_{(2.2.1)}+m_{(2.2.2)}+m_{(2.2.3.1)}
\ge 7m_{(2.1)}+5m_{(2.2.1)}+m_{(2.2.2)}+m_{(2.2.3.1)}
\ge 12m_{(2.2.1)}+8m_{(2.2.2)}+8m_{(2.2.3.1)}
\ge 20m_{(2.2.2)}+8m_{(2.2.3.1)}
\ge 20(m_{(2.2.3.1)}+\cdots+m_{(2.2.3.\ell)})+8m_{(2.2.3.1)}
\ge (20\ell+\nomb8)\almb m_{(2.2.3.\ell)}$, and
$\delta_{(2.2.3.\ell)}=(3\ell+1)m_{(1)}+(2\ell+1)m_{(2.1)}+\ell(m_{(2.2.1)}+m_{(2.2.2)})+m_{(2.2.3.1)}+\cdots+m_{(2.2.3.\ell)}-(8\ell+2)$.
Hence $8\ell+3\le
(3\ell+1)m_{(1)}+(2\ell+1)m_{(2.1)}+\ell(m_{(2.2.1)}+m_{(2.2.2)})+m_{(2.2.3.1)}+\cdots+2m_{(2.2.3.\ell-1)}+m_{(2.2.3.\ell)}
=\frac{3\ell+1}4\mu_{(2.2.3)}-\frac{\ell-1}4m_{(2.1)}+\frac{\ell-1}4(m_{(2.2.1)}+m_{(2.2.2)})-\frac{3(\ell-1)}4m_{(2.2.3.1)}+m_{(2.2.3.2)}+\cdots+2m_{(2.2.3.\ell-1)}+m_{(2.2.3.\ell)}
\le\frac{3\ell+1}4\mu_{(2.2.3)}-(\ell-1)m_{(2.2.3.1)}+m_{(2.2.3.2)}+\cdots+2m_{(2.2.3.\ell-1)}+m_{(2.2.3.\ell)}
\le\frac{3\ell+1}4\mu_{(2.2.3)}+m_{(2.2.3.\ell)}\le(\frac{3\ell+1}4+\frac1{20\ell+8})\mu_{(2.2.3)}$, so
$\mu_{(2.2.3)}\ge(8\ell+3)/(\frac{3\ell+1}4+\frac1{20\ell+8})=4(8\ell+3)(5\ell+2)/(15\ell^2+11\ell+3)\ge\frac{32}3=10.66\cdots$ for any $\ell\ge2$.
On the other hand $\sqrt{321((j-1)k+1)(jk+1)}/9=10.5\cdots$, so (3.30) applies.

(3.39) From now on we may asume $k=2$.
Then $r_{(2.j.2.\ell)}=2j\ell+2j+\ell-2$,
$\mu_{(2.j.2.\ell)}=(j\ell+j-1)m_{(1)}+(\ell+1)(m_{(2.1)}+\cdots+m_{(2.j-1)})+\ell m_{(2/j/1)}+m_{(2.j.2.1)}+\cdots+m_{(2.j.2.\ell)}$ and
$\mu_{(2.j.2)}=(2j-1)m_{(1)}+2(m_{(2.1)}+\cdots+m_{(2.j-1)})+m_{(2.j.1)}+m_{(2.j.2)}
\ge(2j+1)(m_{(2.1)}+\cdots\almb+m_{(2.j-1)})+2jm_{(2.j.1)}+m_{(2.j.2.1)}
\ge((j-1)(2j+1)+2j)m_{(2.j.1)}+((j-1)(2j+1)+1)m_{(2.j.2.1)}
\ge(2j^2+j-1)(m_{(2.j.2.1)}+\cdots+m_{(2.j.2.\ell)})+(2j^2-j)m_{(2.j.2.1)}
\ge(2j-1)((j+1)\ell+j)m_{(2.j.2.\ell)}$.
Hence $2j\ell+2j+\ell-1\le \mu_{(2.j.2.\ell)}+m_{(2.j.2.\ell-1)}
=(j\ell+j-1)m_{(1)}+(\ell+1)(m_{(2.1)}+\cdots\almb+m_{(2.j-1)})+\ell
m_{(2.j.1)}+m_{(2.j.2.1)}+\cdots+2m_{(2.j.2.\ell-1)}+m_{(2.j.2.\ell)}
=\frac{j\ell+j-1}{2j-1}\mu_{(2.j.2)}-\frac{\ell-1}{2j-1}(m_{(2.1)}+\cdots+m_{(2.j-1)})+\frac{(j-1)(\ell-1)}{2j-1}m_{(2.j.1)}-\frac{j(\ell-1)}{2j-1}m_{(2.j.2.1)}+m_{(2.j.2.2)}+\cdots+2m_{(2.j.2.\ell-1)}+m_{(2.j.2.\ell)}
\le\frac{j\ell+j-1}{2j-1}\mu_{(2.j.2)}-(\ell-1)m_{(2.j.2.1)}+m_{(2.j.2.2)}+\cdots+2m_{(2.j.2.\ell-1)}+m_{(2.j.2.\ell)}
\le\frac{j\ell+j-1}{2j-1}\mu_{(2.j.2)}+m_{(2.j.2.\ell)}
\le(\frac{j\ell+j-1}{2j-1}+\frac1{(2j-1)((j+1)\ell+j)})\mu_{(2.j.2)}$, so\nl
$\mu_{(2.j.2)}\ge(2j-1)(2j\ell+2j+\ell-1)((j+1)\ell+j)/(j(j+1)\ell^2+(2j^2-1)\ell+j^2-j+1)$.

(3.40) We divide the case s by $j$.

When $j=8$, $\gamma_{(2.j.2)}=\sqrt{321(2j-1)(2j+1)}/9=31.789\cdots$.
On the other hand
$\mu_{(2.j.2)}\ge15(17\ell+15)(9\ell+8)/(72\ell^2+127\ell+57)\ge\frac{255}8=31.875$ for any $\ell\ge2$ by (3.39), so (3.30) applies.

When $j=7$, we have $\gamma_{(2.j.2)}=27.798\cdots$ and
$\mu_{(2.j.2)}\ge\frac{195}7=27.85\cdots$ for any $\ell\ge2$.
When $j=6$, $\gamma_{(2.j.2)}=23.80\cdots$ and
$\mu_{(2.j.2)}\frac{143}6=23.83\cdots$.
Thus we are done in these cases.

However, the situation is worse when $j\le5$.
For $j=5$, we have $\gamma_{(2.j.2)}=19.807\cdots$, but
$(2j-1)((2j+1)\ell+(2j-1))((j+1)\ell+j)/(j(j+1)\ell^2+(2j^2-1)\ell+j^2-j+1)>\gamma_{(2.j.2)}$ if and only if $\ell\le45$.
Thus the cases $j=5$, $\ell\ge46$ survive.
Similarly, for $j=4$, the cases $\ell\ge5$ survive.
For $j=3, 2$, this method does not work for any $\ell$.

(3.41) We will use (3.23) when $\ell$ is large.
We have $\mu_{(2.j)}=jm_{(1)}+m_{(2.1)}+\cdots+m_{(2.j)}
\ge(j+1)(m_{(2.1)}+\cdots+m_{(2.j-1)}+m_{(2.j.1)})
\ge(j+1)((j-1)m_{(2.j-1)}+m_{(2.j.1)})
\ge(j+1)((j-\nomb1)\almb(m_{(2.j.1)}+m_{(2.j.2.1)})+m_{(2.j.1)})
\ge(j+1)(j(m_{(2.j.2.1)}+\cdots+m_{(2.j.2.\ell)})+(j-1)m_{(2.j.2.1)})
\ge(j+1)(j(\ell-1)+j-1)m_{(2.j.2.\ell-1)}
=(j+1)(j\ell-1)m_{(2.j.2.\ell-1)}$.
Hence $1\le\delta_{(2.j.2.\ell)}+m_{(2.j.2.\ell-1)}$ implies
$2j\ell+2j+\ell-1\le
(j\ell+j-1)m_{(1)}+(\ell+1)(m_{(2.1)}+\cdots+m_{(2.j-1)})+\ell
m_{(2.j.1)}+m_{(2.j.2.1)}+\cdots\almb+m_{(2.j.2.\ell-1)}+m_{(2.j.2.\ell)}
\le(\ell+1-\frac1j)\mu_{(2.j)}+\frac1j(m_{(2.1)}+\cdots+m_{(2.j-1)})-\frac{j-1}j m_{(2.j.1)}+m_{(2.j.2.1)}+\cdots\almb
\le(\ell+1-\frac1j+\frac1{j(j+1)})\mu_{(2.j)}-m_{(2.j.1)}+m_{(2.j.2.1)}+\cdots
\le(\ell+1-\frac1{j+1})\mu_{(2.j)}+m_{(2.j.2.\ell-1)}
\le(\ell+\frac j{j+1}+\frac1{(j+1)(j\ell-1)})\mu_{(2.j)}$, so
$\mu_{(2.j)}\ge((2j+1)\ell+2j-1)/(\ell+\frac j{j+1}+\frac1{(j+1)(j\ell-1)})
=(j+1)(j\ell-\nomb1)((2j+1)\ell+2j-1)/(j\ell+j-1)((j+1)\ell-1)$.

For $j=2$, this implies $\mu_{(2.j)}\ge\sqrt{321j(j+1)}/9$ for $\ell\ge4$, and
(3.23) applies.
This method works when
($j=3$; $\ell\ge4$), ($j=4,5,6$; $\ell\ge3$), ($j\ge7$, $\ell\ge2$).

(3.42) Combining (3.40) and (3.41), we are done except the following cases:\nl
$(j,\ell)=(2,2)$, $(2,3)$, $(3,2)$, $(3,3)$.

In these cases we need a counterpart of (3.23), (3.30) in terms of
$\Xi_{(2.j.2.\ell)}$.
To begin with, let $u$, $w_1$, $w_{2.1}$, $\cdots$ be as before, and set
$\alpha=j\ell+j-1$, $\beta=j\ell+\ell+j$,
$w=(j\ell+j-1)w_1+(\ell+1)(w_{2.1}+\cdots+w_{2.j-1})+\ell
w_{2.j.1}+w_{2.j.2.1}+\cdots+w_{2.j.2.\ell}$ (compare $\mu_{(2.j.2.\ell)}$).
Then\nl
1) $h^0(u,\cdots)=0$ if $\beta u<w$.\nl
2) $h^1(u,\cdots)=0$ if $u\ge w_1+w_{2.1}$, $w_1\ge
w_{2.1}+\cdots+w_{2.j-1}+w_{2.j.1}$, $w_{2.1}\ge \cdots\ge w_{2.j-1}\ge
w_{2.j.1}+w_{2.j.2.1}$, $w_{2.j.2.1}\ge w_{2.j.2.2}\ge\cdots\ge
w_{2.j.2.\ell}$.\nl
3) $h^0(u,\cdots)\le\frac12(u+1)(u+2)-\frac1{2\alpha\beta}w^2+({\text{lower
degree terms in $w$}})$ if $w\le\alpha u$.\nl
4) $h^0(u,\cdots)\le\frac1{2\beta(\beta-\alpha)}(\beta u-w)^2+({\text{lower
degree terms in $\beta u-w$}})$ if $w\le\beta u$.

The proof is similar as before.
Now if $\Xi_{(2.j.2.\ell)}\ge\gamma$, we have
$\Xi_{(2.j.2.\ell)}(\pi_1^*B-xE)\ge\eta_\lambda(x)=\Max(0,\lambda(x-3)+\gamma)$
for some $\lambda\ge\gamma$, so we will estimate
$$\psi(\lambda)=\int_{\eta_\lambda(x)\le\alpha
x}\left(\frac12x^2-\frac{\eta_\lambda(x)^2}{2\alpha\beta}\right)dx+\int_{\alpha
x\le\eta_\lambda(x)\le\beta x}\frac{(\beta
x-\eta_\lambda(x))^2}{2\beta(\beta-\alpha)}dx$$.
As in (3.22) and (3.29), there are three cases:\nl
1) $\frac\gamma3\le\lambda\le3\alpha-\frac23\gamma$.
$6\psi(\lambda)=\frac{729}8-\frac1{\alpha\beta\lambda}(\frac32\lambda+\gamma)^3$.\nl
2) $3\alpha-\frac23\gamma\le\lambda\le3\beta-\frac23\gamma$.
$6\psi(\lambda)=f_{\alpha,\beta,\gamma}(\lambda)$ in (3.11).\nl
3) $3\beta-\frac23\gamma\le\lambda$.
$6\psi(\lambda)=(3\lambda-\gamma)^3/\lambda(\lambda-\alpha)(\lambda-\beta)$.

As before, $\psi(\lambda)$ decreases in the range 1) and 3).
For $c=3\alpha-\frac23\gamma$ and $d=3\beta-\frac23\gamma$, we have by
(3.11)\nl
$\varphi^{'''}(c)=36((5j-8)\ell+5j-13+2\epsilon)>0$,\nl
$\varphi''(d)=50\alpha^2+236(\beta-\alpha)\beta+40(2\beta-\alpha)\epsilon+8\epsilon^2>0$,\nl
$\varphi'(c)=-(\ell+1)(11j\ell-2\ell+11j-13+2\epsilon)(j\ell-\ell+j-2+\epsilon)<0$,\nl
$\varphi'(d)=18(\ell+1)(j\ell+\ell+j)(j\ell+2\ell+j+1-2\ell)>0$,\nl
$\varphi(c)=-36\alpha(\beta-\alpha)(2\alpha-\beta+\epsilon)^2<0$ and\nl
$\varphi(d)=-9\beta(\beta-\alpha)(4\alpha^2+7\beta(\beta-\alpha)+8(2\beta-\alpha)\epsilon+4\epsilon^2)<0$\nl
unless $\epsilon$ is very large.
{}From this we infer that $\varphi<0$ in the range 2).
Hence
$\psi(\lambda)\le\psi(\frac\gamma3)=\frac{27}{16}(9-\frac{\gamma^2}{\alpha\beta})$.
Thus we conclude

(3.43) {\sl $\pi_1^*B-\frac92E$ is big if
$\Xi_{(2.j.2.\ell)}\ge\sqrt{321\alpha\beta}/9$, where $\alpha=j\ell+j-1$ and
$\beta=j\ell+\ell+j$.}

(3.44) We have $\mu_{(2.j.2.\ell)}=\alpha
m_{(1)}+(\ell+1)(m_{(2.1)}+\cdots+m_{(2.j-1)})+\ell
m_{(2.j.1)}+m_{(2.j.2.1)}+\cdots+m_{(2.j.2.\ell)}
\ge(\alpha+\ell+1)(m_{(2.1)}+\cdots+m_{(2.j-1)})+(\alpha+\ell)m_{(2.j.1)}+m_{(2.j.2.1)}+\cdots+m_{(2.j.2.\ell)}
\ge\beta(j-\nomb1)m_{(2.j-1)}+(\beta-1)m_{(2.j.1)}+m_{(2.j.2.1)}+\cdots+m_{(2.j.2.\ell)}
\ge(\beta j-1)m_{(2.j.1)}+(\beta j-\beta+1)
m_{(2.j.2.1)}\almb+m_{(2.j.2.2)}+\cdots+m_{(2.j.2.\ell)}
\ge((\beta j-1)(\ell-1)+\beta j-\beta+\ell-1)m_{(2.j.2.\ell-1)}
=\beta(\alpha-j)m_{(2.j.2.\ell-1)}$.
Hence $1\le\delta_{(2.j.2.\ell)}+m_{(2.j.2.\ell-1)}$ implies
$\mu_{(2.j.2.\ell)}
\ge(\alpha+\beta)/(1+\frac1{\beta(\alpha-j)}=(\alpha+b)\beta(\alpha-j)/\alpha(\beta-j-1)$.
Therefore (3.43) applies when $\ell=3$ and $j=2,3$.
But the case $\ell=2$ still survive.

(3.45) Similarly as in (3.37), we infer that the blow-up process continues via
$E_{(2.j.2.2.n)}$ for some $n\ge2$.
$$(1)-(2.j.1)-(2.j.2.2.1)-\cdots-(2.j.2.2.n)-(2.j.2.1)-(2.j-1)-\cdots-(2.1)$$
We have
$\mu_{(2.j.2)}=(2j-1)m_{(1)}+2(m_{(2.1)}+\cdots+m_{(2.j-1)})+m_{(2.j.1)}+m_{(2.j.2.1)}
\ge(2j+1)(m_{(2.1)}+\cdots+m_{(2.j-1)})+2jm_{(2.j.1)}+m_{(2.j.2.1)}
\ge((j-1)(2j+1)+2j)m_{(2.j.1)}+((j-1)(2j+1)+1)m_{(2.j.2.1)}
\ge(2j+1)(2j-1)m_{(2.j.2.1)}+((j-1)(2j+1)+2j)m_{(2.j.2.2.1)}
\ge((2j+1)(2j-1)(n-1)+2j^2+j-1)m_{(2.j.2.2.n-1)}=(2j-1)((2j+1)n-j)m_{(2.j.2.2.n-1)}$.
Hence $1\le\delta_{(2.j.2.2.n)}+m_{(2.j.2.2.n-1)}$ implies $4jn+2j+1
\le((2j-1)n+j)m_{(1)}+(2n+1)(m_{(2.1)}+\cdots+m_{(2.j-1)})+(n+1)m_{(2.j.1)}+nm_{(2.j.2.1)}+m_{(2.j.2.2.n-1)}+\cdots+2m_{(2.j.2.2.n-1)}+m_{(2.j.2.2.n)}
\le(n+\frac
j{2j-1})\mu_{(2.j.2)}-\frac1{2j-1}(m_{(2.1)}+\cdots+m_{(2.j-1)})+\frac{j-1}{2j-1}m_{(2.j.1)}-\frac j{2j-1}m_{(2.j.2.1)}+m_{(2.j.2.2.1)}+\cdots
\le(n+\frac j{2j-1})\mu_{(2.j.2)}-m_{(2.j.2.1)}+m_{(2.j.2.2.1)}+\cdots
\le(n+\frac j{2j-1})\mu_{(2.j.2)}+m_{(2.j.2.2.n-1)}
\le(n+\frac j{2j-1}+\frac1{(2j-1)((2j+1)n-j)})\mu_{(2.j.2)}$,
so $\mu_{(2.j.2)}\ge(2j-1)(4jn+2j+1)((2j+1)n-j)/((2j+1)n+j+1)((2j-1)n-j+1)$.

For $j=3$, $\mu\ge5(12n+7)(7n-3)/(7n+4)(5n-2)=_{\text{def}}\varphi(n)$.
Since $\frac d{dn}\log\varphi=(49n^2+126n+22)/(12n+7)(7n-3)(7n+4)(5n-2)>0$,
we have $\varphi(n)\ge\varphi(2)=\frac{1705}{144}=11.84\cdots$.
On the other hand $\sqrt{321(2j-1)(2j+1)}/9=11.77\cdots$, so this case is ruled
out by (3.30).

For $j=2$, we have $\mu\ge3(8n+5)(5n-2)/(5n+3)(3n-1)=_{\text{def}}\varphi(n)$.
Similarly $\frac d{dn}\log\varphi>0$ and
$\varphi(n)\ge\varphi(2)=\frac{504}{65}=7.75\cdots$.
On the other hand $\sqrt{321(2j-1)(2j+1)}/9=7.71\cdots$, thus this case is
ruled out either.
\dnl
{\bf \S4. Concluding remarks}

(4.1) Combining the observations in preceding sections, we can now complete the
proof of the following

\Th {\sl Let $B$ be a nef and big line bundle on a smooth threefold $M$.
Let $K$ be the canonical bundle and $x$ be a point on $M$.
Suppose that $BC\ge3$ for any curve $C\ni x$, $B^2S\ge7$ for any surface $S\ni
x$ and $B^3\ge51$.
Then $K+B$ is spanned at $x$.}

In the proof, we should pay attention to the danger (3.32).
We argue as follows.

For the sake of simplicity, at first, we make the assumption below:
$$\Bs\vert s(\pi_1^*B-\epsilon E)\vert\cap E=\emptyset\text{ for some
$\epsilon>0$ and $s\gg 0$.}\tag*$$
In such a case $\Xi_P(\pi_1^*B-\epsilon E)=0$ at any place $P$ of $E$, so the
assertions in (3.23), (3.30) etc. are valid even if
$\Xi_P=\sqrt{321\alpha\beta}/9$, since $\lambda>\frac\gamma3$.

Assume that $x\in\Bs\vert K+B\vert$.
Then $\pi_1^*B-\frac92E$ is not big by Theorem 0.
Take a Hironaka decomposition as in (1.4).
There are only finitely many non-exceptional bad places.
By the results in \S2, we can replace the Hironaka model by another one having
no non-exceptional bad places.

Let us call the union of the loci of bad places as ``bad locus'', and a point
on it as a ``bad point''.
This notion may depend on the degree of the Hironaka model, but it does not
depend on the choice of the model $\tilde\pi:\tilde M\lra M$.
Moreover, if we replace the model by another one whose degree is a multiple of
the older one, then the new bad locus is a subset of the older one.
In particular, no new bad place appear by such a replacement.

There are at most finitely many bad points $p_1,\cdots,p_n$ on $E\cong\BP^2$.
By (3.10;1), we can choose a Hironaka model of degree $\tau$ such that
$\mu_{(1)}=\tau^{-1}v_{p_j}(\tau\vert\pi_1^*B-3E\vert)<1.99072$ at each point
$p_j$.
Then the primary place $\SP(Y_{(1)})$ over $p_j$ is not bad, so the bad locus
on the $(-1)$-curve $Y_{(1)}$ over $p_j$ is a finite set.
By (3.10;2), we can choose a model such that $\mu_{(2)}<2.81531$ at any such
points on $Y_{(1)}$ over $p_j$.
Let $Y_{(2)}$ be the $(-1)$-curve over a bad point on $Y_{(1)}$.
Then $\SP(Y_{(2)})$ is not a bad place, so there are only finitely many bad
points on it.
By (3.10;3), we may assume $\mu_{(3)}<3.51$ at every such point.
Similarly, there appear only finitely many bad points by such a blowing-up, so
we may assume $\mu_{(4)}<4.23$ and $\mu_{(5)}<5$ by (3.10).

If $\tilde Y_{(2)}\cap Y_{(3)}$ is a bad point, we apply (3.19) and assume
$\mu_{(3.2)}<6.31$, $\mu_{(3.3)}<9.13$ and $\mu_{(3.4)}<12$ at every bad point
of this type.
Similarly, if $\tilde Y_{(1)}\cap Y_{(2)}$ is a bad point, we arrange things so
that\nl
$\mu_{(2.j)}<\sqrt{321j(j+1)}/9$ for any $j\le108$ (cf. (3.23)),\nl
 $\mu_{(2.j.k)}<\sqrt{321((j-1)k+1)(jk+1)}/9$ for any $j\le8$, $k\le12$ (cf.
(3.30)),\nl
$\mu_{(2.j.2.\ell)}<\sqrt{321(j\ell+j-1)(j\ell+\ell+j)}/9$ for any $j\le3$,
$\ell\le3$.\nl
Up to this step, there are only finitely many places involved.

Now, by the argument in \S3, we infer that there are no bad places on such a
model.
Hence (1.6) applies.

Even if the hypothesis $(*)$ is not true, we may assume that the numbers
$\mu_{(2.j)}$, $\mu_{(2.j.k)}$ and $\mu_{(2.j.2.\ell)}$ are bounded from above
by numbers just slightly bigger than $\sqrt{321\alpha\beta}/9$, which are
enough for the later argument in \S3.
Thus the same method works to complete the proof.

(4.2) \Cor {\sl $K+3L$ is spanned if $L$ is an ample line bundle with $L^3>1$.}

(4.3) By the same method we can prove the following

\Th {\sl Let $(M,\Delta)$ be a log threefold having only log terminal
singularities.
Let $A$ be a line bundle on $M$ such that $B=A-K$ is a log ample $\BQ$-bundle,
where $K=K(M,\Delta)=K_M+\Delta$ is the log canonical $\BQ$-bundle of
$(M,\Delta)$.
Let $x$ be a smooth point on $M$ off $\Delta$ such that $BC\ge3$ for any curve
$C\ni x$, $B^2S\ge7$ for any surface $S\ni x$ and $B^3\ge51$.
Then $A$ is spanned at $x$.}

Probably the last assumption can be replaced by a weaker one $B^3\ge d$, as
long as $d>50.625$, if one works harder than in this paper.
Any way, we can hardly think that these numbers $51$ or $50.625$ are optimal.
For, $\pi_1^*B-3E$ is big if $B^3>27$, so $x\in\Bs\vert L\vert$ implies in this
case that $\pi_1^*B-3E$ fails to be nef.
Hence there exists some troublesome subvariety containing $x$.
A better understanding of the nature of it will yield a better numerical
criterion of spannedness.
A true optimist may hope that it is enough to assume $BC\ge3$, $B^2S\ge9$ and
$B^3>27$.

Probably there is a similar numerical criterion even if $x$ is a singular point
on $(M,\Delta)$.

(4.4) The method in this paper will be used in a forthcoming paper to study the
very ampleness of $K+B$.
I hope that the idea may have some value in higher dimensions too.

\enddocument